\RequirePackage{fix-cm}
\documentclass[twocolumn]{svjour3}          % twocolumn
\smartqed  % flush right qed marks, e.g. at end of proof
\usepackage{graphicx}
 \usepackage{mathptmx}      % use Times fonts if available on your TeX system
\usepackage{bm}
\usepackage{fancyhdr}
\usepackage{stmaryrd}
\usepackage{amsmath}
\usepackage{comment}
\usepackage{natbib}
\usepackage{float}
\usepackage{stfloats}
\usepackage{amssymb}
\usepackage{dashrule}

\usepackage{latexsym}

\SetSymbolFont{letters}{bold}{OML}{cmm}{b}{it}
\SetSymbolFont{operators}{bold}{OT1}{cmr}{bx}{n}

 \DeclareMathAlphabet{\mathcal}{OMS}{cmsy}{m}{n}
\newcommand{\vect}[1]{\boldsymbol{\mathbf{#1}}}
\begin{document}\sloppy

\title{A penalized approach to the bivariate logistic regression model for the association between ordinal responses
}


\author{Marco Enea         \and
        Gianfranco Lovison %etc.
}


\institute{Marco Enea \and Gianfranco Lovison \at
              Dipartimento di Scienze Economiche Aziendali e Statistiche, University of Palermo, Palermo, Italy; \\
              Tel.: +39-09123895334\\
              Fax: +39-091485726\\
              \email{marco.enea@unipa.it; gianfranco.lovison@unipa.it}
}

\date{Received: date / Accepted: date}

\maketitle

\begin{abstract}
Bivariate ordered logistic models (BOLMs) are appealing to jointly model the marginal distribution of two ordered responses and their association, given a set of covariates. When the number of categories of the responses increases, the number of global odds ratios (or their re-parametrizations) to be estimated also increases and estimating the association structure becomes crucial for this type of data. In fact, such data could be too \lq\lq rich" to be fully modelled with an ordinary BOLM while, sometimes, the well-known Dale's model could be too parsimonious to provide a good fit.  In addition, when the cross-tabulation of the responses contains some zeros, for a number of model configurations, including the bivariate version of the partial proportional odds model (PPOM), estimation of a BOLM by the Fisher-scoring algorithm may either fail or estimate a too \lq\lq irregular" association structure. \par In this work, we propose to use a nonparametric approach for the maximum likelihood estimation of a BOLM. We apply penalties to the differences between adjacent row and column effects. As a result, estimation is less demanding than an ordinary BOLM, permitting the fit of PPOMs and/or the smoothing of the marginal and association parameters by polynomial curves and surfaces, with scores chosen by the data. Model selection is based on the penalized log-likelihood ratio, whose limiting distribution has been studied through simulations, and AIC. Our proposal is compared to the Goodman's model and the Dale's model, in terms of goodness-of-fit and parsimony, on a literature data set. Finally, an application on an original data set of liver disease patients is proposed.
\keywords{Dale model \and bivariate ordered logistic model \and penalized maximum likelihood estimation \and ordinal association}
\end{abstract}

\section{Introduction}
\label{intro}
\noindent Models for association play a central role in ordered categorical data analysis. For the multivariate case, marginal models (MMs) represent a natural choice to model marginal distributions of the responses given covariates. An example of full likelihood based marginal model is \cite{Dale:1986}. A similar model, the multivariate logistic model described in \citet{Gloneck:1995}, but restricted to the bivariate ordered version, is the basis on which we develop our proposal. Some open, or at least not completely solved, problems about estimation of a multivariate ordered logistic model are of computational type and concern  maximum likelihood (ML) estimation by iterative algorithms, often providing invalid estimates at the $k$th step, exceeding the boundaries of the parameter space. Some of such problems could be solved as in \citet{Colombi:2001} and \citet{Bartolucci:2002} by including strict inequality constraints. However, constrained ML estimation is appealing only when a particular application implies natural ordering constraints. On the contrary, when the ordering is not fully reliable, or externally imposed, like in responses which arise from discretized versions of latent continuous variables, using inequality constraints may not be appropriate. Indeed, due to lack of subject-matter knowledge that yields natural restrictions on marginal distributions, no strict ordering constraints are appropriate, and more helpful and flexible approaches are necessary.
In these situations, a nonparametric approach may be useful \citep{Dardanoni:1998}. Within the possible range of nonparametric approaches, penalization is the one considered in this paper. Surprisingly, there is little literature on penalization applied to marginal models.  \citet{Desantis:2008} apply a ridge penalty to a latent class model for ordinal data to stabilize ML estimation, that would otherwise not be computationally feasible without application of strict constraints. Other contributions deal mainly with forms of longitudinal  \citep{Gieiger:1997,Fahrmeier:1999} or horizontal \citep{Bustami:2001} nonparametric modeling. The former focuses on smoothing of variation of marginal and association parameters over time, the latter refers to a form of smoothing on covariates, often by using splines.\par Our proposal is based on a form of vertical smoothing - that is across response levels - of the regression parameters in order to regularize the parameter space and/or fit polynomial models using scores ``chosen by the data''. After recalling the Dale, the Gloneck-McCullagh and the bivariate partial proportional odds models in Section 2, and the penalized ML estimation approach in Section 3, the penalty terms we propose to use are introduced in Section 3. By simulation, we show the advantage of using the proposed approach in Section 4 and check the asymptotic behaviour of the penalized deviance statistic in Section 5. Two applications are considered in Section 6: in the first one, we compare our proposal to the \citet{Dale:1986} and the \citet{Goodman:1979} models on a literature data set, whereas the second one is about a data set of liver disease patients.

\section{Bivariate ordered logit models}\label{sec2}
\noindent For two ordered outcomes  $A_1$ and $A_2$, define the row and column marginal cumulative probabilities of a $D_1 \times D_2$ contingency table $A_1A_2$ as
\begin{eqnarray}
\mu_{r.}=P(A_{1}\leq r)=\sum_{i \leq r}\pi_{i.},\,\,\,\,\,
\mu_{.c}=P(A_{2}\leq c)=\sum_{j \leq c}\pi_{.j},\nonumber
\label{margprob}
\end{eqnarray}
and the upper-left quadrant probabilities as
\begin{eqnarray}
\mu_{rc}=P(A_{1}\leq r,\,A_{2}\leq c)&=&\sum_{i\leq r}\sum_{j\leq c}\pi_{ij},\nonumber
\label{quadrantprob}
\end{eqnarray}
with $r=1,\dots,D_{1},\,c=1,\dots,D_{2}$. By differencing we obtain
\begin{eqnarray}
P(A_{1}\leq r,\,A_{2}> c)&=&\mu_{r.}-\mu_{rc},\nonumber\\
P(A_{1}> r,\,A_{2}\leq c)&=&\mu_{.c}-\mu_{rc},\nonumber\\
P(A_{1}> r,\,A_{2}> c)&=&1-\mu_{r.}-\mu_{.c}+\mu_{rc}.\nonumber
\end{eqnarray}
By choosing the cumulative odds as ordinal risk measures, and the logit as link function, we obtain the \emph{global logits} (or \emph{log global odds}):
\begin{eqnarray}
\log\phi_{1r}&=& \textrm{logit} [P(A_{1}\leq r)]=\log(\mu_{r.})-\log(1-\mu_{r.}),
\label{l1}
\end{eqnarray}
\begin{eqnarray}
\log\phi_{2c}&=& \textrm{logit} [P(A_{2}\leq c)]=\log(\mu_{.c})-\log(1-\mu_{.c}),
\label{l2}
\end{eqnarray}
$r=1,\dots,D_{1}-1,\,c=1,\dots,D_{2}-1$. By choosing the cross-products of quadrant probabilities as ordinal association measures, and the natural logarithm as link function, the \emph{log global odds ratios} (or log-GORs) are defined as:\\
\begin{eqnarray}
\log\psi_{rc}&=&\log\frac{P(A_{1}\leq r,\,A_{2}\leq c)P(A_{1}> r,\,A_{2}> c)}{P(A_{1}\leq r,\,A_{2}> c)P(A_{1}> r,\,A_{2}\leq c)}\nonumber\\
&=&\log\frac{\mu_{rc}(1-\mu_{r.}-\mu_{.c}+\mu_{rc})}
       {(\mu_{r.}-\mu_{rc}) (\mu_{.c}-\mu_{rc})}.
\label{lor}
\end{eqnarray}
Given the three parameters $\mu_{r.}$, $\mu_{.c}$, and $\psi_{rc}$, we may find the corresponding joint cumulative probabilities with the following inversion formula:\\
\begin{eqnarray}
\mu_{rc}=
\begin{cases}
\frac{1}{2}(\psi_{rc}-1)^{-1}(a_{rc}-\sqrt{a_{rc}^2+b_{rc}})&  \textrm{if}\, \psi\ne 1,\\
\mu_{r.}\mu_{.c}& \textrm{if}\, \psi=1,
\end{cases}
\label{plack}
\end{eqnarray}
where $a=1+(\mu_{r.}+\mu_{.c})(\psi_{rc}-1)$ and $b=-4\psi_{rc}(\psi_{rc}-1)\mu_{r.}\mu_{.c}$.
If the cumulative probabilities $\mu_{r.}$ and $\mu_{.c}$ satisfy the constraints $\mu_{r.}<\mu_{r+1,.}$ for $r=1,\dots,D_{1}-1$, and $\mu_{.c}<\mu_{.,c+1}$ for $c=1,\dots,D_{2}-1$, and the global odds ratios are not dependent on the category, that is $\psi_{rc}=\psi$, then (\ref{plack}) is a Plackett distribution \citep{Plackett:1965}. Thus, the bivariate Dale regression model for $(\bm{\phi}_1,\bm{\phi}_2,\bm{\psi}_{12})'$ is as follows:
\begin{equation}
\begin{cases}
\log[\phi_{1,r}(\bm{x})]=\beta_{10r}-\bm{\beta}_{1}'\bm{x},\\
%&\\
\log[\phi_{2,c}(\bm{x})]=\beta_{20c}-\bm{\beta}_{2}'\bm{x},\\
%&\\
\log[\psi_{rc}(\bm{x})]=\alpha+\rho_{1r}+\rho_{2c}+\sigma_{rc}-\bm{\beta}_{3}'\bm{x},\\
\end{cases}
\label{Dalemodel}
\end{equation}
$r=1,...,D_{1}-1,\,\,c=1,...,D_{2}-1$. This model does not require marginal scores for responses and it is also invariant under any monotonic transformation of the marginal responses. Further, since the model is based on global odds ratios, collapsing adjacent row or column categories does not produce any effect in parameter interpretation, which remains unchanged with the exception of the intercepts related to the collapsed categories. This is in contrast with the RC Goodman model which uses local cross-ratios. In a more general framework than (\ref{Dalemodel}), \citet{Gloneck:1995} introduce the \emph{multivariate logistic model}:
\begin{equation}
\vect{C}'\log(\vect{L}\bm{\pi})=\vect{X}\bm{\beta},
\label{mod2}
\end{equation}
where $\vect{C}$ is a contrasts matrix, $\vect{L}$ is a matrix with elements $a_{ij}\ge 0$ such that $\vect{L}\bm{\pi}=\bm{\mu}$, $\bm{\eta}=\vect{C}'\log(\vect{L}\bm{\pi})$ is the parameter vector of interest,  and $\vect{X}$, an $n \times p$ matrix, with $n=\prod_{k=1}^{K}{D_{k}}$. Although formulation (\ref{mod2}) is referred to $K\geq2$ responses, here only two responses $A_{1}$ and $A_{2}$ are considered.  The components of  $\vect{C}'\log(\vect{L}\bm{\pi})$ are symbolically denoted by
$\bm{\eta}=$$(\eta_{\varnothing},\bm{\eta}_{A_{1}}',\bm{\eta}_{A_{2}}',\bm{\eta}_{A_{1}A_{2}}')'$,
where $\eta_{\varnothing}=\log(\sum \bm{\pi})=0$ is the null contrast and the remaining vectors have  elements specified by (\ref{l1}),  (\ref{l2}) and (\ref{lor}), respectively. We will refer to (\ref{mod2}) as the bivariate ordered logistic model (BOLM).  \citet{Lapp:1998} show how to fit the Dale and Goodman models starting from the framework of a BOLM.
\noindent Some computational problems may arise when fitting a multivariate logistic model, depending on the number of responses and categories.  For example, when inverting equation $\bm{\eta}=\vect{C}'{\log}(\vect{L}\bm{\pi})$ to obtain $\bm{\pi}$ in terms of $\bm{\eta}$, it may happen that for certain fixed values of $\bm{\eta}$ no positive solution $\bm{\pi}$ exists. Although $\bm{\pi}>0$ ensures the matrix $\vect{C}'\vect{D}^{-1}\vect{L}$ to be invertible \citep[Theorem 1]{Gloneck:1995}, where $\vect{D}=\mbox{diag}(\vect{L}\bm{\pi})$, the range of the mapping is not a hyper-rectangle and fixing some components of $\bm{\eta}$ restricts the range of the remaining components, that is the model is not \emph{variation independent}. Although this problem is particularly magnified for $K>2$ responses \citep{Bergsma:2002a,Qaqish:2006},
computational problems can also arise in the bivariate case, above all when considering certain particular model configurations. For instance, it may happen that not only the intercepts but some covariates have a category-dependent effect. To highlight this effect one may want to fit a bivariate version of the partial proportional odds model proposed by \citet{Peterson:1990}. However, such a model can be computationally very hard to fit, even with a limited and reasonable number of parameters. To deal with this difficulty, we propose to regularize the parameter space by penalizing the log-likelihood of the model. This allows to increase the range of possible models to be fitted. The penalty term we use for this is introduced in Section \ref{section:3.3.2}.
\par The fit of a BOLM becomes computationally hard also when the number of response categories increases. In addition, the model may result overparameterized. \citet{Lapp:1998} fit a Dale's model by imposing constraints on the row and column interactions of the association intercepts in order to reduce the number of parameters. However, this type of data appears to be too ``rich" to be modeled with fully parametric models and nonparametric or semiparametric models, followed by graphical presentation, could result more useful \citep{Eilers:1996}. In order to smooth the marginal and association effects across the response categories, we suggest to use a penalty term, introduced is Section \ref{section:3.3.3}, for nonparametric modeling, mainly employed in the P-spline context \citep{Eilers:2006}, but suitably re-written to be used in the framework of a BOLM. In part, this approach can be considered the bivariate extension of the models proposed by \cite{Tutz:2003b}.
\par The ordinal nature of the responses imposes inequality constraints on marginal distributions which have to be taken into account in model estimation. In Section \ref{section:3.3.1}, we present a penalty term, which is able to mimic such inequality constraints.

%\subsection{The non-proportional, the proportional and the partially proportional odds model}\label{section:3.1}
\noindent In order to better understand the potential of the penalization approach, some further notation is needed,
according to that used in \citet{Tutz:2003} for the univariate cumulative logistic regression model. Let $\mathcal{Q}$ be the set of indices of all the   covariates, excluding the intercepts, and $\mathcal{P}\subset \mathcal{Q}$ be a subset of $p$ covariates. Let $\mathcal{S}$ be the set of indices of the variables whose effects we assume do not depend on categories and such that $\mathcal{S}\subseteq\mathcal{P}$, and let $\bar{\mathcal{S}}=\mathcal{P} \backslash \mathcal{S}$. In particular, we define $\mathcal{S}$ and $\bar{\mathcal{S}}$ as $\mathcal{S}=\cup_{k=1}^{3}\mathcal{S}_{k}$, and $\bar{\mathcal{S}}=\cup_{k=1}^{3}\bar{\mathcal{S}}_{k}$, where $\mathcal{S}_{k}$ and $\bar{\mathcal{S}}_{k}$ are the subsets of $\mathcal{S}$ and $\bar{\mathcal{S}}$, respectively, associated to the $k$th equation. To complete the notation, let $\mathcal{S}^{0}=\{0\}\cup\mathcal{S}$ and $\bar{\mathcal{S}}^{0}=\{0\}\cup\bar{\mathcal{S}}$. Consider the following model where only a part of the covariates is supposed to be category-independent:

\begin{equation}
\label{sis2}
\begin{cases}
\log[\phi_{1r}(\bm{x}_{i})]=\beta_{10r}+\bm{\beta}_{1\mathcal{S}_{1}}'\bm{x}_{i\mathcal{S}_{1}}+
\bm{\beta}_{1\bar{\mathcal{S}}_{1}r}'\bm{x}_{i\bar{\mathcal{S}}_{1}},\\
%&\\
\log[\phi_{2c}(\bm{x}_{i})]=\beta_{20c}+\bm{\beta}_{2\mathcal{S}_{2}}'\bm{x}_{i\mathcal{S}_{2}}+
\bm{\beta}_{2\bar{\mathcal{S}}_{2}c}'\bm{x}_{i\bar{\mathcal{S}}_{2}},\\
%&\\
\log[\psi_{rc}(\bm{x}_{i})]=\beta_{30rc}+\bm{\beta}_{3\mathcal{S}_{3}}'\bm{x}_{i\mathcal{S}_{3}}+
\bm{\beta}_{3\bar{\mathcal{S}}_{3}rc}'\bm{x}_{i\bar{\mathcal{S}}_{3}},\\
\end{cases}
\end{equation}
$(r=1,\ldots,D_{1}-1,\,\,c=1,\ldots,D_{2}-1)$. We refer to model (\ref{sis2}) as the \emph{Non-Uniform association and Partially Proportional Odds Model} (NUPPOM). Although in the univariate case the phrase ``proportional odds'' is usually referred to a model with covariate effects which do not depend on the categories, here we will refer to a \emph{Uniform association and Proportional Odds Model} (UPOM) as a model defined from (\ref{sis2}) assuming $\beta_{30rc}=\beta_{30}$ and $\bar{\mathcal{S}}=\varnothing$. On the other hand, a \emph{Non-Uniform association and Non-Proportional Odds Model} (NUNPOM) will be defined from (\ref{sis2}) assuming $\mathcal{S}=\varnothing$ and with category-dependent association intercepts. Note that the intercepts for the marginal equations (that is the global-logit intercepts) are never supposed to be independent of the categories, whatever the model. According to these definitions the bivariate Dale model (\ref{Dalemodel}) is a NUPOM, and it becomes a UPOM when $\rho_{1r}=0$, $\rho_{2c}=0$ and $\sigma_{rc}=0$, $r=1,\ldots,D_{1}-1,\,\,c=1,\ldots,D_{2}-1$. Further, to specify that a NUPPOM is fitted we will also write $NUPPOM(\mathcal{S})$ and to indicate that a UPOM is fitted we will also write $UPOM(\bar{\mathcal{S}}^0)$.

\indent Under multinomial sampling with frequencies $\bm{y}_{i}\sim M(n_{i},$ $\bm{\pi}_{i})$, consider the model  $\vect{C}'\log(\vect{L}\bm{\pi}_{i})=\vect{X}_{i}\bm{\beta}$, with the matrices $\vect{C}$ and $\vect{L}$ such that the marginal parameters are \emph{global logits}, the association parameters are \emph{log global odds ratios}, and the constraint $\sum_{j=1}^{D_{1}}\sum_{k=1}^{D_{2}}{\pi_{ijk}=1}$ is included. Then, the kernel of the log-likelihood is
\begin{equation}
l(\bm{\beta})=\sum_{i=1}^{m}{l(\bm{\beta};\bm{y}_{i})} = \sum_{i=1}^{m}{\bm{y}_{i}'\log(\bm{\pi}_{i})},
\end{equation}
where $m$, the observed number of response configurations, is such that $\sum_{i=1}^{m}{n_i}=n$, with $n$ indicating the sample size. The \emph{penalized log-likelihood} has the form

\begin{equation}
l_{P}(\bm{\beta})=l(\bm{\beta})-\frac{1}{2}\tau(\bm{\beta}),
\label{lp}
\end{equation}

\noindent where $\tau(\bm{\beta})=\bm{\beta}'\vect{P}\bm{\beta}$, and $\vect{P}$ represents the penalization and includes the smoothing parameter. Penalized ML estimation formulas are given and discussed in Appendix A. The specification of $\tau(\bm{\beta})$ is discussed in the following sections, whereas the form of $\vect{P}$ is given in Appendix B.

\subsection{A penalty term for parameter space regularization}\label{section:3.3.2}
\noindent When the cross-tabulation of the responses contains one or more zeros, parameter estimation by Fisher-scoring may be challenging at each iteration. In these cases, one may try to reduce $l_{step}$, the step length (see Appendix A). However, estimates of the association structure may result too irregular, with very high (or very low) estimated odds ratios in correspondence of the zero cells. In order to stabilize the ML estimates of the BOLM using Fisher scoring, a reduction of the parameter space may be helpful. We propose to penalize both the marginal and the association parameters. In addition, since the model is not variation independent, applying a penalty term on association parameters might be useful to limit the range of the possible values that the marginal parameters can assume, so avoiding a failure of the Fisher scoring. The general expression is: \\
\begin{eqnarray}
\tau(\bm{\beta})&=&\sum_{j\in \bar{\mathcal{S}}_{1}^0}\lambda_{1j}\sum_{r=2}^{D_1-1}{\zeta( \beta_{1jr})}+\sum_{j\in\bar{\mathcal{S}}_{2}^0}\lambda_{2j}\sum_{c=2}^{D_2-1}{
\zeta(\beta_{2jc})}\nonumber\\
&+&\sum_{j\in\bar{\mathcal{S}}_{3}^0}\lambda_{3j}\sum_{r=2}^{D_1-1}
\sum_{c=2}^{D_2-1}{\zeta(\beta_{3jrc})},
\label{sis3}
\end{eqnarray}
\noindent where $\lambda_{kj}$ is the smoothing parameter for the $j$th variable of the $k$th equation of system (\ref{sis2}), $k=1,2,3$. A first specification of $\zeta(.)$, that we call ARC1, is  $\zeta(\alpha_t)=(\Delta\alpha_{t})^2$, where $\Delta$ is the order 1 difference operator, that is $\Delta\alpha_{t}=\alpha_{t}-\alpha_{t-1},\,t\geq 2$. With respect to (\ref{sis3}), the operator acts over the indices $r$ and/or $c$. It involves penalization of adjacent row and column parameters and is aimed at ($i$)  overcoming estimation problems by reducing parameter space and ($ii$) reproducing a $UPOM$ for high smoothing values. For the aim $(i)$, the choice of $\lambda$ is based on the minimum value for which Fisher-scoring does not fail. The simulation study in Section \ref{section:3.5} will clarify this choice. The aim ($ii$) is achieved as $\lambda_{kj}\shortrightarrow\infty,\,k=1,2,3,\,\forall j \in \bar{\mathcal{S}}_{k}^0$, for which all the parameters indexed by $j$ will tend to be equal among the categories.  Although (\ref{sis3}) allows to penalize marginal intercepts, it is preferable to avoid a strong penalization on such parameters, in order not to violate (\ref{constr}).
In this work emphasis is on ARC1, but several other specifications of $\zeta(.)$ are possible. For example, another specification aimed at reducing the parameter space  is $\zeta(\alpha)=\alpha^2$, which corresponds to a ridge-type penalty for the bivariate logistic regression model. As $\lambda_{kj}\shortrightarrow\infty,\,k=1,2,3,\,\forall j \in \bar{\mathcal{S}}_{k}^0$, all the parameters indexed by $j$ will tend to zero. A similar penalty, not involving the third term, is used by \cite{Desantis:2008} in a penalized latent class model for ordinal data.

\subsection{A penalty term for nonparametric modeling}\label{section:3.3.3}
\noindent Beside being useful for reducing the parameter space and for reproducing a $UPOM$, the following generalization of the penalty term ARC1, hereafter denoted by ARC2, can be used to specify row or column effects and to fit nonparametric models where the effects are determined by a polynomial:

%\begin{multline}
\begin{eqnarray}
\tau(\bm{\beta})\hspace{-.1cm}&=&\hspace{-.3cm}\sum_{j\in\bar{\mathcal{S}}_{1}^0}\lambda_{1\hspace{-.05cm}j}
\sum_{r=s_{1}+1}^{D_1-1}{\hspace{-.1cm}(\Delta^{s_{1j}}\beta_{1jr})^2}+
\hspace{-.1cm}\sum_{j\in\bar{\mathcal{S}}_{2}^0}\lambda_{2\hspace{-.03cm}j}\sum_{c=s_{2}+1}^{D_2-1}{\hspace{-.1cm}(\Delta^{s_{2j}} \beta_{2jc})^2} \nonumber\\
\hspace{-.3cm}&+&\hspace{-.3cm}\sum_{j\in\bar{\mathcal{S}}_{3}^0}\left[ \lambda_{3j}\hspace{-.3cm}\sum_{r=s_{3}+1}^{D_1-1}\sum_{c=1}^{D_2-1}{(\Delta^{s_{3j}} \beta_{3jrc})^2}\right. \nonumber\\
&+&\left.\lambda_{4j}\hspace{-.1cm}\sum_{r=1}^{D_1-1}\sum_{c=s_{4}+1}^{D_2-1}{(\Delta^{s_{4j}} \beta_{3jrc})^2}\right],\hspace{.3cm}
\label{sis4}
\end{eqnarray}
%\end{multline}
\noindent where $\Delta^{a}=\Delta(\Delta^{a-1})$. Consider the following penalty settings:

\begin{itemize}
\item as $\lambda_{hj}=0, \,h=1,...,4,\,\forall j \in \bar{\mathcal{S}}^0$, an unrestricted model will be fitted;

\item as $\lambda_{hj}\shortrightarrow\infty, \,h=1,...,4,\,\forall j \in \bar{\mathcal{S}}_{1}\cup\bar{\mathcal{S}}_{2}\cup\bar{\mathcal{S}}_{3}^0$ and $s_{hj}=1$, the fitted parameters will tend to be equal, and the model will tend to a $UPOM$;

\item  as $\lambda_{3j}\shortrightarrow\infty$, $\lambda_{4j}=0,\,\forall j \in \bar{\mathcal{S}}_{3}^0$ and $s_{3j}=1$, a model with column effects will be fitted;

\item as $\lambda_{3j}=0$, $\lambda_{4j}\shortrightarrow\infty,\,\forall j \in \bar{\mathcal{S}}_{3}^0$ and $s_{4j}=1$, a model with row effects will be fitted;

\item  as $\lambda_{hj}\shortrightarrow\infty, \,h=1,2,\,\forall j \in \bar{\mathcal{S}}^0$ and $s_{hj}>1$, the fitted parameters will follow a polynomial curve of degree $s_{hj}-1$.
\item  as $\lambda_{hj}\shortrightarrow\infty, \,h=3,4,\,\forall j \in \bar{\mathcal{S}}^0$ and $s_{hj}>1$, the fitted parameters will follow a polynomial surface of degree $s_{3j}+s_{4j}-2$.
\end{itemize}
%In particular, for the latter point, while all the estimated marginal parameters, indexed by $j$, will lie onto a polynomial curve, all the estimated association parameters, indexed by $j$, will lie onto a polynomial surface.
Notice the difference between the penalty terms included in (\ref{sis4}) and those included in the penalized log-likelihood (14) in \cite{Tutz:2003b}, suggested for a single ordered response. In that paper, the author proposed to penalize the differences of adjacent categories, for a vertical smoothing, jointly to the use of penalized B-splines for a horizontal smoothing, resulting in a form similar to (\ref{sis4}). Also notice the differences with the bivariate horizontal smoothing approach by \cite{Bustami:2001} which presented the additive bivariate Dale model, for continuous, category-independent covariates, as a natural extension of the generalized additive model \citep{Hastie:1990}.
\par Penalty (\ref{sis4}) may be useful to assume certain dependence structures on the categories, for both marginal and association parameters. For example, if one wants to assume a linear trend for the row marginal effects, one may assume $\eta_{1ir}=\beta_{10r}+\sum_{j=1}^{p}x_{ij}\beta_{1jr}$, where $\beta_{1jr}=\alpha_{0j}+\alpha_{1j}\delta_{jr}$, with $\alpha_{0j}$ and $\alpha_{1j}$ unknown parameters, and with scores $\delta_{jr}$. In spite of its simplicity, such an approach assumes arbitrary scores. An alternative way is just to use a penalization approach with penalty term ARC2 which uses scores ``chosen by the data'' \citep{Tutz:2003}. Indeed, the smoothing parameters and the polynomial degrees can be chosen on the basis of some criterion, such as the values that minimize the AIC. As a special case, suppose to want to fit a model which assumes a linear trend of the marginal parameters and an association structure composed by the interaction of two first degree polynomials\footnote{The degree of a two-variable polynomial is defined as the highest degree of its terms, and the degree of a term is the sum of the exponents of the variables that appear in it. Since (\ref{sis4}) allows to fit only polynomial models with interactions, to distinguish each of the possible models having the same degree, it is more practical for us to indicate a model by specifying both the degrees of the one-variable polynomials, omitting to specify the (implicit) presence of interaction terms.}. By assuming, for simplicity, that the same variable $x_{j}$ is present in all the equations of system (\ref{sis2}), choosing $s_{hj}=2, \,h=1,\ldots,4$ and high smoothing values, for instance $10^8$, the predictor becomes
\begin{itemize}
\item $x_{ij}\beta_{jr}=x_{ij}\gamma_{01j}+x_{ij}(\gamma_{11j}\hspace{-0.05cm}\cdot\hspace{-0.05cm} r)$
\item $x_{ij}\beta_{jc}=x_{j}\gamma_{02j}+x_{ij}(\gamma_{12j}\hspace{-0.05cm}\cdot\hspace{-0.05cm} c)$
\item $x_{ij}\beta_{jrc}=x_{ij}\gamma_{03j}+x_{j}(\gamma_{13j}\hspace{-0.05cm}\cdot\hspace{-0.05cm} r)+
x_{ij}(\gamma_{23j}\hspace{-0.05cm}\cdot\hspace{-0.05cm} c)+x_{ij}(\gamma_{33j}\hspace{-0.05cm}\cdot r\hspace{-0.05cm}\cdot\hspace{-0.05cm} c)$,
\end{itemize}

\noindent with scores $\delta_{jr}=r,\,\delta_{jc}=c,\,$ and $\bm{\delta}_{jrc}=(r,c)'$, that is pre-assigned equally-spaced scores. %Because of the interaction term, the third point represent a second degree polynomial, but here we will indicate this model to be linear for row and column.

%or over the predictor for the $i-th$ observation
%\begin{equation}
%\tau(\bm{\beta})=\sum_{i=1}^{n}\left[\lambda_{1}\sum_{r=1}^{d1-1}{ %|\bm{\eta}_{i,r}|}+\lambda_{2}\sum_{c=1}^{d2-1}{|\bm{\eta}_{i,c}|}+
%\lambda_{3}\sum_{r=1}^{d1-1}\sum_{c=1}^{d2-1}{|\bm{\eta}_{i,rc}|}\right].
%\label{sis_lasso2}
%\end{equation}
\subsection{Mimicking inequality constraints}\label{section:3.3.1}
\noindent
The ordinal nature of the responses introduces some explicit ordering constraints on marginal distribution which have to be taken into account to avoid ill-conditioning of the predictor space. In particular, for the $i$th individual, such  constraints are on the marginal predictors, that is
\begin{equation}
\beta_{k01}\hspace{-0.05cm}+\hspace{-0.05cm}\bm{\beta}_{k1}'\bm{x}_{i}\hspace{-0.05cm}<\hspace{-0.05cm}
\beta_{k02}\hspace{-0.05cm}+\hspace{-0.05cm}\bm{\beta}_{k2}'\bm{x}_{i}\hspace{-0.05cm}<\hspace{-0.05cm}
...\hspace{-0.05cm}<\hspace{-0.05cm}\beta_{k0,D_{k}-1}\hspace{-0.05cm}+\hspace{-0.05cm}
\bm{\beta}_{k,D_{k}-1}'\bm{x}_{i},
%\alpha_{21}+\bm{\beta}_{1}'\bm{x}_{i}<\alpha_{22}+\bm{\beta}_{2}'\bm{x}_{i}<...<\alpha_{2c}+
%\bm{\beta}_{c}'\bm{x}_{i}.
\label{constr}
\end{equation}
$k=1,2.$ Although Lagrangians can be used to take into account such constraints, in the spirit of this paper, a penalized-oriented solution could be the following:

%\begin{equation}
%\tau(\bm{\beta})=\sum_{i=1}^{n}\left[\lambda_{1}\sum_{r=2}^{D1-1}{ I(\Delta\eta_{1ir})(\Delta\eta_{1ir})^2}+\lambda_{2}\sum_{c=2}^{D2-1}{I(\Delta\eta_{2ic})
%(\Delta\eta_{2ic})^2}\right],
%\label{sis_constraints}
%\end{equation}
\begin{equation}
\tau(\bm{\beta})=\sum_{i=1}^{n}\left[\sum_{k=1}^{2}\lambda_{k}\sum_{r=2}^{D_k-1}{ I(\Delta\eta_{kir})(\Delta\eta_{kir})^2}\right],
\label{sis_constraints}
\end{equation}

\noindent where $\eta_{kir}=\beta_{k0r}+\bm{\beta}_{kr}'\bm{x}_{i}$, $\Delta\eta_{kir}=\eta_{kir}-\eta_{ki,r-1}$, and
$I(z)=1$ if $z\geq 0$, otherwise $I(z)=0$.  As $\lambda_{k}\shortrightarrow\infty,$ the penalty term (\ref{sis_constraints}) acts in such a way to satisfy (\ref{constr}).The univariate version of (\ref{sis_constraints}) is used, for example, by \citet{Muggeo:2008} in a penalized splines context applied to univariate generalized linear models. It can also be used jointly to (\ref{sis3}) or (\ref{sis4}). Notice that, although seemingly superfluous, the inclusion of $(\Delta\eta_{kir})^2$  in (\ref{sis_constraints}) derives from the necessity of writing $\tau(\bm{\beta})$ as a quadratic form in order to exploit the penalized ML formulae in Appendix A.

\section{Hypothesis testing}\label{section:testing}
\noindent When estimates are penalized, the asymptotic distribution of the penalized likelihood ratio ($LR_p$) statistic is known only for some hypothesis systems. Unfortunately, to check the hypothesis $(P)POM$ to which we are mainly interested, i.e. that of category-independent effects, as far as we know, neither exact nor asymptotic results are known for $LR_P$. Thus, by a simulation study, we analyze the conditions under which it is possible to approximate, under the hypothesis ($P)POM$, the $LR_p$ asymptotic distribution by a $\chi^2$ distribution. As an introduction, in the following Section \ref{gray} we first recall a result already present  in the literature, which is useful for a simple hypothesis system, and then we show the rationale of using the $\chi^2$ distribution and the results of the simulation study in Section \ref{simlrp}. To simplify notation we will suppose, without loss of generality, that the same index $j$ refers to the same variable for both marginal and association equations.

\subsection{The $LR_P$ statistic for the hypothesis of null effects\label{gray}}
\noindent Let us consider the specific partition of parameters  $\bm{\beta}_{\mathcal{P}^0}=(\bm{\gamma},\bm{\delta})'$, such that the null hypothesis:
\begin{equation}
H_0: \bm{\delta}=\bm{0},
\label{HS1}
\end{equation}
postulates that only a subset of parameters is constrained. Furthermore, consider the penalized log-likelihood of the more general model, $l_{P}(\hat{\bm{\gamma}},\hat{\bm{\delta}})$, that of the reduced model, $l_{P}(\tilde{\bm{\gamma}},\bm{0})$ and the penalized log-likelihood ratio statistic:
\begin{equation}
LR_{P}=-2\{l_{P}(\tilde{\bm{\gamma}},\bm{0})-l_{P}(\hat{\bm{\gamma}},\hat{\bm{\delta}})\}.
\end{equation}
\noindent Let $\mathbf{F}$ be the information matrix from the unpenalized partial likelihood, with subscripts denoting the submatrices, such as $\mathbf{F}_{\bm{\delta}\bm{\delta}}$ for derivatives with respect to $\bm{\delta}$. Consider the matrix $\mathbf{F}_{\bm{\delta}\bm{\delta}|\bm{\gamma}}=\mathbf{F}_{\bm{\delta}\bm{\delta}}-
\mathbf{F}_{\bm{\delta}\bm{\gamma}}\mathbf{F}^{-1}_{\bm{\gamma}\bm{\gamma}}\mathbf{F}_{\bm{\gamma}\bm{\delta}}$. Then, under the null hypothesis, \cite{Gray:1994} shows the statistic $LR_{P}$ to have the same asymptotic distribution as
$\sum{\alpha_jZ^2_j}$, where the $Z_j$'s are independent standard Normal random variables, and the $\alpha_j$'s are the eigenvalues of the matrix $\lim_{n \to \infty}\mathbf{F}_{\bm{\delta}\bm{\delta}|\bm{\gamma}}(\mathbf{F}_{\bm{\delta}\bm{\delta}|\bm{\gamma}}+\mathbf{P})^{-1}$, where $\mathbf{P}$ is the matrix representing the penalty term.

\subsection{The $LR_P$ statistic for the $(P)POM$ hypothesis\label{simlrp}}
Consider a full model of the $NUNPOM$ type, i.e. for which all variables $j$, $j \in \mathcal{P}^0\equiv\{\bar{\mathcal{S}}^0,\mathcal{S}=\varnothing\}$, have category-dependent effects, and a reduced model for which the effects of some variables $j$, $j \in \mathcal{S}\ne\varnothing$, are category independent.
The penalized log-likelihood ratio test to check the hypothesis for comparing these two models, i.e. for testing the null hypothesis:
\begin{equation}
H_{0}: \bm{\beta}_j=\beta_j\bm{1},\,\,j \in \mathcal{S},
\label{HS2}
\end{equation}
compares the maximum penalized log-likelihood $l_{P}(\hat{\bm{\beta}}_{\mathcal{P}^0})$, and the maximum penalized log-likelihood $l_{P}(\tilde{\bm{\beta}}_{\mathcal{S}},\tilde{\bm{\beta}}_{\bar{\mathcal{S}}^0})$:

\begin{eqnarray}
LR_{P}&=&-2\{l_{P}(\tilde{\bm{\beta}}_{\mathcal{S}},\tilde{\bm{\beta}}_{\bar{\mathcal{S}}^0})-
l_{P}(\hat{\bm{\beta}}_{\mathcal{P}^0})\}\nonumber\\
&=&2\sum_{i=1}^m\sum_{r=1}^{D_{1}}\sum_{c=1}^{D_{2}} \bm{y}'_{irc}\log\left(\frac{\hat{\bm{\pi}}_{irc}}{\tilde{{\bm{\pi}}}_{irc}}\right)+\tau(\tilde{\bm{\beta}})-
\tau(\hat{\bm{\beta}}),
\end{eqnarray}

\noindent where $\hat{\bm{\pi}}'_{i}=(\hat{\pi}_{i11},\ldots,\hat{\pi}_{iD_{1}D_{2}})'$ is the estimated (by penalization) probability vector for the model under $H_1$ and $\tilde{\bm{\pi}}_{i}=(\tilde{\pi}_{i11},\ldots,\tilde{\pi}_{iD_{1}D_{2}})'$ is the corresponding estimated (by penalization) probability vector for the reduced model. %As far as we know, neither exact nor asymptotic results are known for this statistic under this hypothesis system. However, under certain conditions, we might approximate the asymptotic distribution by a $\chi^2$ with appropriate degrees of freedom. We show this
Supposing to use the penalty term ARC1 and following \citet{Tutz:2003}, let $\lambda_{kj\mathcal{R}}\,(\lambda_{kj\mathcal{F}})$ denote the smoothing parameters for the reduced model (full model). Then, we have
\begin{displaymath}
\tau(\tilde{\bm{\beta}})-\tau(\hat{\bm{\beta}})=
\end{displaymath}
\begin{eqnarray}
&=&\sum_{j\in \bar{\mathcal{S}}^0}\left[\lambda_{1j\mathcal{R}}\sum_{r=2}^{D_1-1}{(\Delta\tilde{\beta}_{1jr})^2}+
\lambda_{2j\mathcal{R}}\sum_{c=2}^{D_2-1}{(\Delta\tilde{\beta}_{2jc})^2}\right.\nonumber\\
&+&\left.\lambda_{3j\mathcal{R}}\sum_{r=2}^{D_1-1}
\sum_{c=2}^{D_2-1}{(\Delta \tilde{\beta}_{3jrc})^2}\right]-\sum_{j\in\mathcal{P}^0}\left[\lambda_{1j\mathcal{F}}\sum_{r=2}^{D_1-1}{(\Delta\hat{\beta}_{1jr})^2}\right.\nonumber\\
&+&\left.\lambda_{2j\mathcal{F}}
\sum_{c=2}^{D_2-1}{(\Delta\hat{\beta}_{2jc})^2}+\lambda_{3j\mathcal{F}}\sum_{r=2}^{D_1-1}\sum_{c=2}^{D_2-1}{(\Delta \hat{\beta}_{3jrc})^2}\right]\nonumber\\
&=&\sum_{j\in\bar{\mathcal{S}}^0}\left\{\sum_{r=2}^{D_1-1}\left[\lambda_{1j\mathcal{R}}(\Delta\tilde{\beta}_{1jr})^2-
\lambda_{1j\mathcal{F}}(\Delta\hat{\beta}_{1jr})^2\right]\right.\nonumber\\
&+&\left.\sum_{c=2}^{D_2-1}\left[\lambda_{2j\mathcal{R}}
(\Delta\tilde{\beta}_{2jc})^2-\lambda_{2j\mathcal{F}}(\Delta\hat{\beta}_{2jc})^2\right]\right.\nonumber\\
&+&\left.\sum_{r=2}^{D_1-1}\sum_{c=2}^{D_2-1}
\left[\lambda_{3j\mathcal{R}}(\Delta \tilde{\beta}_{3jrc})^2-\lambda_{3j\mathcal{F}}(\Delta \hat{\beta}_{3jrc})^2\right]\right\}\nonumber\\
&-&\sum_{j\in\mathcal{S}}\left[\lambda_{1j\mathcal{F}}\sum_{r=2}^{D_1-1}{
(\Delta\hat{\beta}_{kjr})^2}+\lambda_{2j\mathcal{F}}\sum_{c=2}^{D_2-1}{
(\Delta\hat{\beta}_{kjc})^2}\right.\nonumber\\
&+&\left.\lambda_{3j\mathcal{F}}\sum_{r=2}^{D_1-1}
\sum_{c=2}^{D_2-1}{(\Delta \hat{\beta}_{kjrc})^2}\right].\nonumber
\end{eqnarray}

\noindent If estimates are not penalized, that is if for $k=1,2,3$, $\lambda_{k1\mathcal{R}}=\lambda_{k2\mathcal{R}}=\dots=\lambda_{k|\bar{\mathcal{S}}^0|\mathcal{R}}=\lambda_{k0\mathcal{F}}=
\lambda_{k1\mathcal{F}}=\dots=\lambda_{k|\mathcal{P}^0|\mathcal{F}}=0$,
 one obtains $\nu(\tilde{\bm{\beta}})-\nu(\hat{\bm{\beta}})=0$, and the $LR_{P}$ statistic has the usual asymptotic $\chi^2$ distribution. If $\lambda_{kj\mathcal{R}}=\lambda_{kj\mathcal{F}}$ is chosen for $j\in\bar{\mathcal{S}}^0$,   $k=1,2 ,3$, then the first term is very small since $\tilde{\beta}_{kjr}\approx\hat{\beta}_{kjr},\,\tilde{\beta}_{kjc}\approx\hat{\beta}_{kjc}$ and $\tilde{\beta}_{kjrc}\approx\hat{\beta}_{kjrc}$ for $r=1,\dots,D_{1}-1,\,c=1,\dots,D_{2}-1$. Thus, the fundamental term concerns the variables for which $j \in \mathcal{S}$ and, if estimates are penalized with a low smoothing value, converging to zero at an appropriate rate, the same asymptotic behaviour holds. \par We show this approximate result by simulation. We simulated the sampling distribution of $LR_P$ assuming to have responses $A_{1}$ and $A_{2}$, with levels $D_{1}=D_{2}=3$, and a single dichotomous covariate $X_{1}$, sampled from $Ber(0.5)$. The chosen regression parameters for the true model are the following:
\begin{itemize}
\item $\bm{\beta}_{10}=(-0.5,0.5)',\,\,\bm{\beta}_{11}=(-0.3,0.3)'$,
\item $\bm{\beta}_{20}=(-0.1,0.6)',\,\,\bm{\beta}_{21}=(-0.2,0.4)'$,
\item $\bm{\beta}_{30}=(1.5,2,2.5,3)',\,\,\bm{\beta}_{31}=(-.5,-.5,-.5,-.5)'$.
\end{itemize}

\begin{figure*}[t]
 \begin{center}
 % \includegraphics[width=15cm,height=9cm]{LRp.pdf}
 %\vspace{-.7cm}
  \includegraphics[angle=-90, scale=.55]{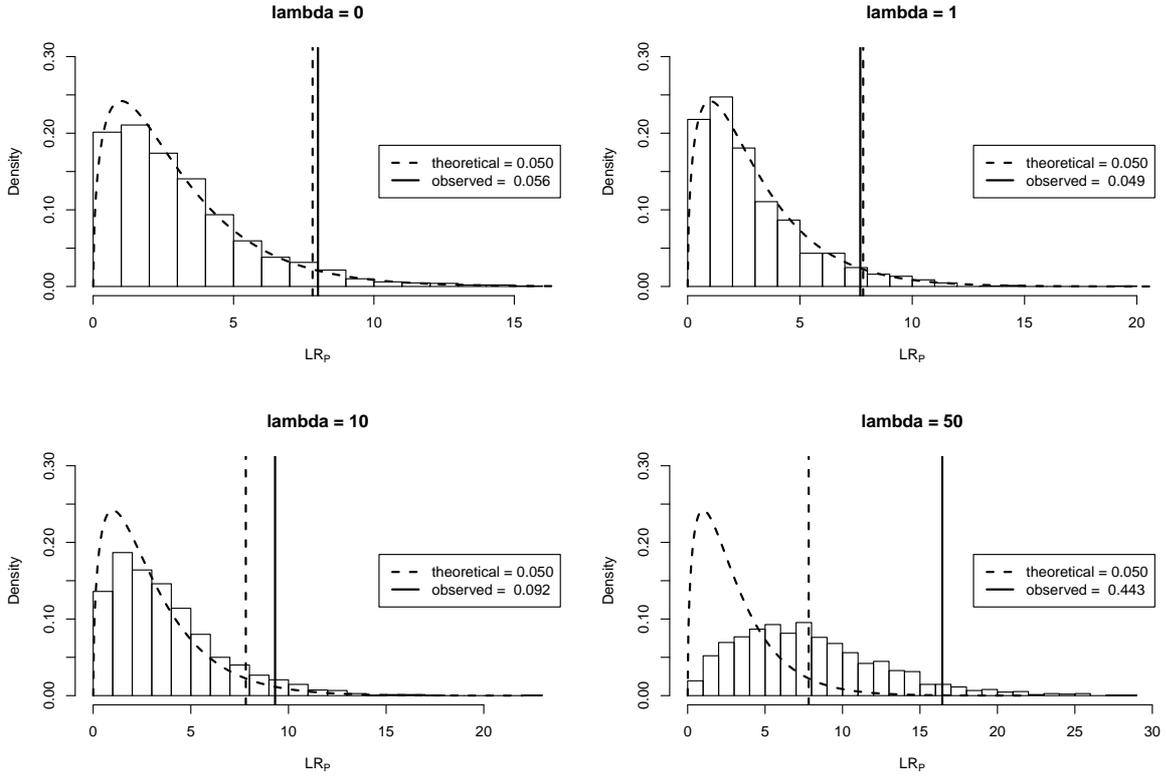}
 \caption{Simulated distribution of $LR_P$ for varying smoothing parameter and a sample size of 400 simulated observations.}\label{trv400}
  \end{center}
\end{figure*}

\noindent  The null hypothesis is that of global association effect for variable $X_1$, that is $H_{\mathcal{S}_{3}}:\beta_{311}=\beta_{312}=\beta_{313}=\beta_{314}=-.5$. We generated 1500 pseudo-samples of size $n=400$ from a multinomial distribution with probability matrix $\boldsymbol{\mathrm{\Pi}}(\mathbf{X}\bm{\beta})$ and fitted the reduced and unreduced models by penalizing the association intercepts $\bm{\beta}_{30}$ according to smoothing values $\lambda_{30} \in \{0,1,10,50\}$. The sampling distribution of $LR_P$, for varying smoothing parameters, is displayed in Figure \ref{trv400} together with the superimposed theoretical $\chi_3^2$ distribution, and two vertical lines highlighting the nominal 5\% level and the observed level. Observe that the theoretical $\chi_3^2$ approximation to the empirical  $LR_{P}$ distribution is very good for $\lambda_{30}\leq 1$, whereas it is totally wrong for $\lambda_{30}=50$.

\section{Evaluating the performance of penalized estimates}\label{section:3.5}
\noindent In order to evaluate the potential of smoothed estimates, a small simulation study is carried out. The sample size is taken to be $n=400$, whereas the number of samples is taken to be $N=100$. A $NUNPOM$ is assumed to hold, with responses $A_{1}$ and $A_{2}$ with $D_{1}=D_{2}=3$ levels, and with
\begin{itemize}
\item $\bm{\beta}_{10}=(-0.6,0.6)',\,\,\bm{\beta}_{11}=(0.3,-0.3)'$,
\item $\bm{\beta}_{20}=(-0.6,0.6)',\,\,\bm{\beta}_{21}=(-0.6,0.6)'$,
\item $\bm{\beta}_{30}=(2.6,2.4,2.0,1.7)',\,\,\bm{\beta}_{31}=(-0.4,0.2,-0.5,0.5)'$.
\end{itemize}

\noindent The covariate values $x_i,\,i=1,...,n$ were drawn from a uniform distribution on $(-1,1)$. Thus, given the model formula (\ref{mod2}) and the inversion method (\ref{plack}) we found the $n\times (D_{1}D_{2})$ probability matrix $\boldsymbol{\mathrm{\Pi}}$, in which each row $\bm{\pi}'_{i}$ represents the probability vector for the $i$th observation. Then, each of $n$ bivariate responses was drawn from a multinomial distribution with probability vector $\bm{\pi}'_{i}$. The $UPOM$ was compared
to the $NUNPOM$, for which ARC1 has been used in combination with (\ref{sis_constraints}). Penalization parameters for ARC1, that is $\lambda_{1},\lambda_{2}$ and $\lambda_{3}$ were chosen equal to a single value $\lambda$, varying in the set $\{0,1,10,100,...\}$. By starting from the lower $\lambda$ value, the procedure was stopped when Fisher scoring did not fail. This procedure was iterated for each of $N$ samples. Comparisons were made by evaluating the following loss functions:\\

\par Mean squared error loss:
\begin{equation}
MSEL=\frac{1}{n}\sum_{i=1}^{n}\sum_{r=1}^{D_{1}D_{2}}(\pi_{ir}-\hat{\pi}_{ir})^2,
\end{equation}
\par Mean relative squared error loss:
\begin{equation}
MRSEL=\frac{1}{n}\sum_{i=1}^{n}\sum_{r=1}^{D_{1}D_{2}}\frac{(\pi_{ir}-\hat{\pi}_{ir})^2}{\pi_{ir}},
\end{equation}
\par Mean entropy or Kullback-Leibler loss:
\begin{equation}
MEL=\frac{1}{n}\sum_{i=1}^{n}\sum_{r=1}^{D_{1}D_{2}}\pi_{ir}\textrm{log}\left(\frac{\pi_{ir}}{\hat{\pi}_{ir}}\right),
\end{equation}
and the AIC defined in Appendix A. Fisher scoring without penalization failed in 45 out of 100 simulations when the $NUNPOM$ was assumed, that is when $\lambda=0$. Before setting $\lambda$ to some value greater than zero, some attempts to estimate the model were made by reducing the step length $l_{step}$ of the iterative algorithm (see Appendix A), and in some case the $NUNPOM$ was fitted. The $UPOM$ was fitted in all the simulations. The results are reported in Table \ref{tabsim1}.

\begin{table}[h]
\small
\centering
\caption{Comparison $UPOM$ vs $NUNPOM$ with ARC1, in terms of simulated mean values of the loss functions and $AIC$. $^*$Cumulative number of Fisher scoring successes ($FSS$).}
%\vspace{.2cm}
\begin{tabular}{rrrrrrr}
{\it Model} & {\it $\lambda$} & {\it MSEL} & {\it MRSEL} &  {\it MEL} &  {\it AIC} & {\it FSS$^*$} \\
\hline
NUNPOM &  0 &   0.0042 &     0.0408 &     0.0208 &     1592.8 &         55 \\

NUNPOM &  1 &   0.0048 &     0.0415 &     0.0225 &     1579.0 &         68 \\

NUNPOM &  10 &   0.0046 &     0.0448 &     0.0227 &     1572.2 &         93 \\

NUNPOM &  100 &   0.0093 &     0.0927 &     0.0399 &     1582.9 &       100 \\

 UPOM &   -  &  0.0129 &     0.1330 &     0.0531 &     1608.8 &        100 \\
\hline
\end{tabular}
\label{tabsim1}
\end{table}
%\vspace{-.4cm}

\noindent Observe that for the $NUNPOM$ all the mean loss functions and $AIC$ are smaller than $UPOM$ ones. The $NUNPOM$ without penalization (i.e. $\lambda=0$) has the smallest loss functions values, as it should be, but also the greatest $AIC$ value among the $NUNPOMs$. This is due to the trade-off between the best fitting, given by the saturated model, and the greater flexibility, given by the penalization approach, which reduces the degrees of freedom.

\section{Applications to real data sets}\label{section:5}
\subsection{The British males occupational status data set}\label{section:3.6}
\noindent Consider the data on occupational status (OS) of a sample of British males from \citet{Goodman:1979}, where fathers and their sons were cross-classified according to the occupational status using seven ordered categories. The data are reported in Table \ref{socialstatus}.
% Table generated by Excel2LaTeX from sheet 'Foglio1'
\begin{table}[h]
\small
\caption{Cross-classification of British males according to occupational status.}
\vspace{.1cm}
\centering
\begin{tabular}{lrrrrrrr}
\hline
  Father's &            &             \multicolumn{4}{r}{Subject's status}             &            &            \\

    status &            &            &            &            &            &            &            \\

           &          1 &          2 &          3 &          4 &          5 &          6 &          7 \\
\hline
         1 &         50 &         19 &         26 &          8 &         18 &          6 &          2 \\

         2 &         16 &         40 &         34 &         18 &         31 &          8 &          3 \\

         3 &         12 &         35 &         65 &         66 &        123 &         23 &         21 \\

         4 &         11 &         20 &         58 &        110 &        223 &         64 &         32 \\

         5 &         14 &         36 &        114 &        185 &        714 &        258 &        189 \\

         6 &          0 &          6 &         19 &         40 &        179 &        143 &         71 \\

         7 &          0 &          3 &         14 &         32 &        141 &         91 &        106 \\
\hline

\end{tabular}
%\vspace{.5cm}
\label{socialstatus}
\end{table}

%\vspace{-.2cm}
\noindent Several authors have re-analyzed such data. For example, \citet{Lapp:1998} compare the Goodman RC and Dale models in terms of goodness-of-fit. We further re-analyze the data by fitting the BOLM with ARC2. The aim of the application is to show the advantages of our proposal when compared to the existing alternatives.

\vspace{-.7cm}
\begin{figure}[H]
 \begin{center}
  \includegraphics[width=8cm,height=7cm]{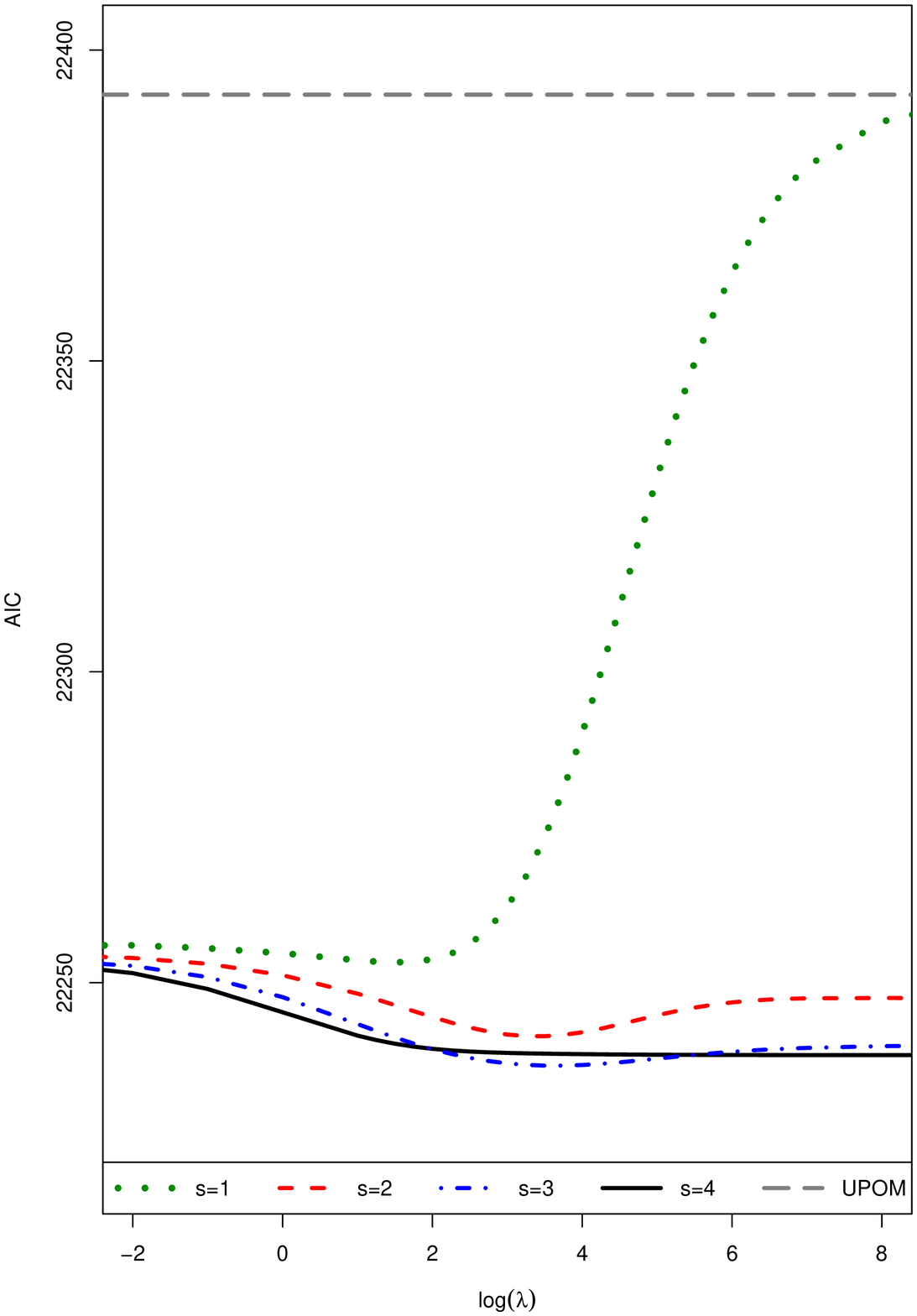}\\
  \includegraphics[width=8cm,height=7cm]{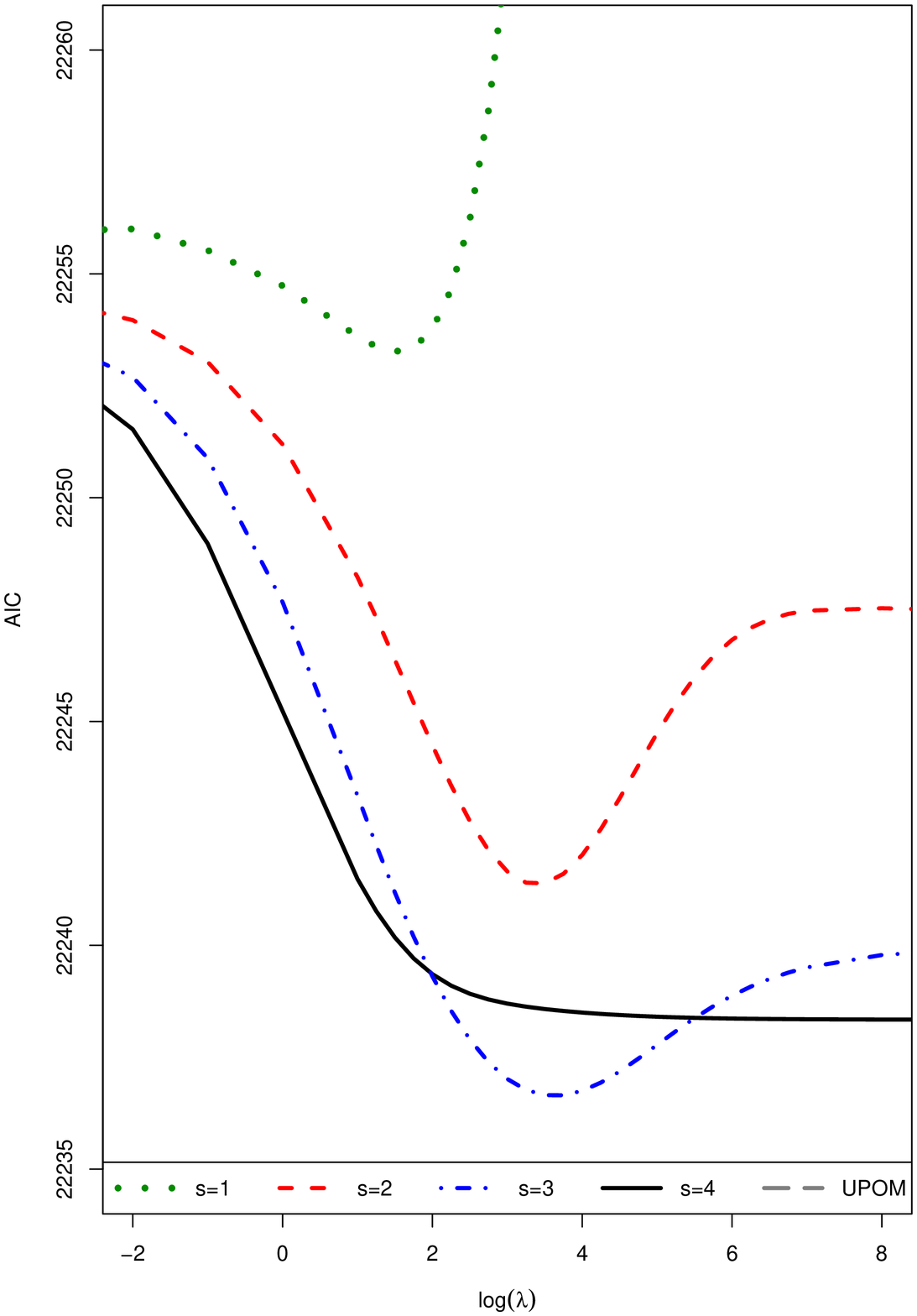}\\
  % \vspace{-.3cm}
  \caption{OS data set, the top graph: AIC for the NUPOM using the penalty term ARC2 for varying smoothing parameter and different orders of penalization. The bottom graph shows a detail of the most critical interval (-2,8).}\label{aic4}
  \end{center}
\end{figure}
\vspace{-.7cm}

\begin{figure*}[t]
\centering%
\vspace{-1cm}
\includegraphics[width=7cm, height=7.5cm]{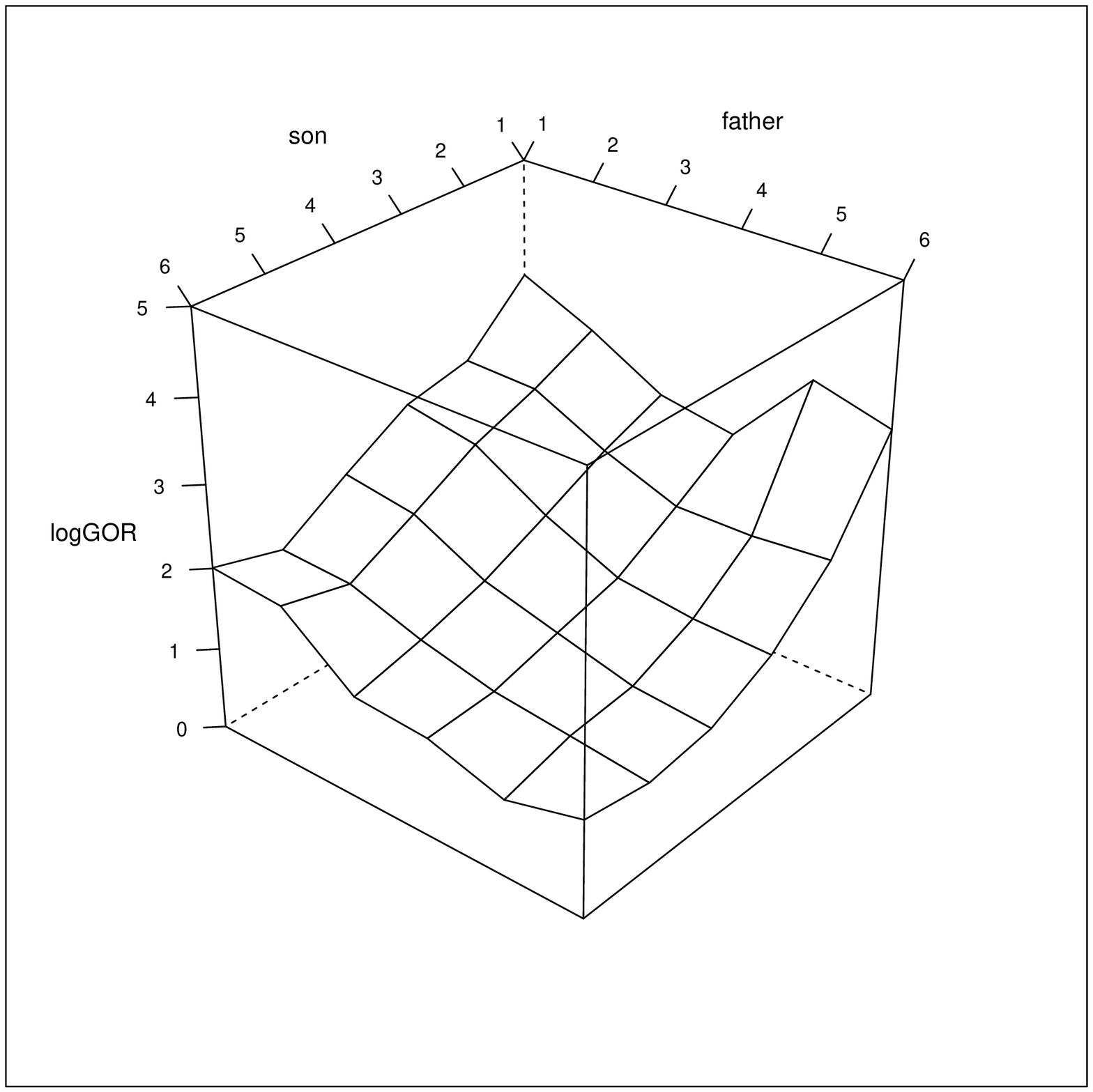}\hspace{-.15cm}
\includegraphics[width=7cm, height=7.5cm]{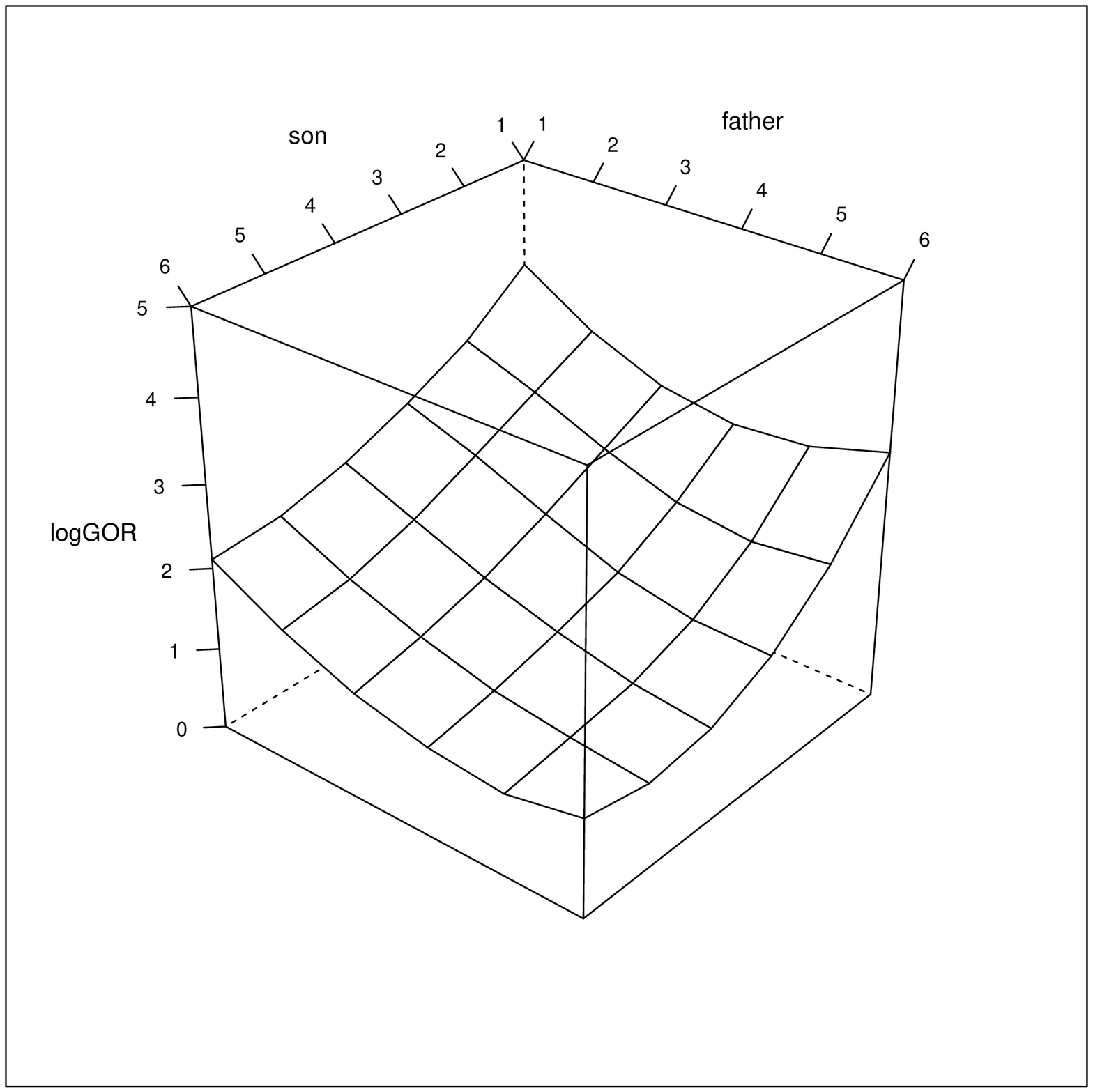}\vspace{-2.35cm}\\
\includegraphics[width=7cm, height=7.5cm]{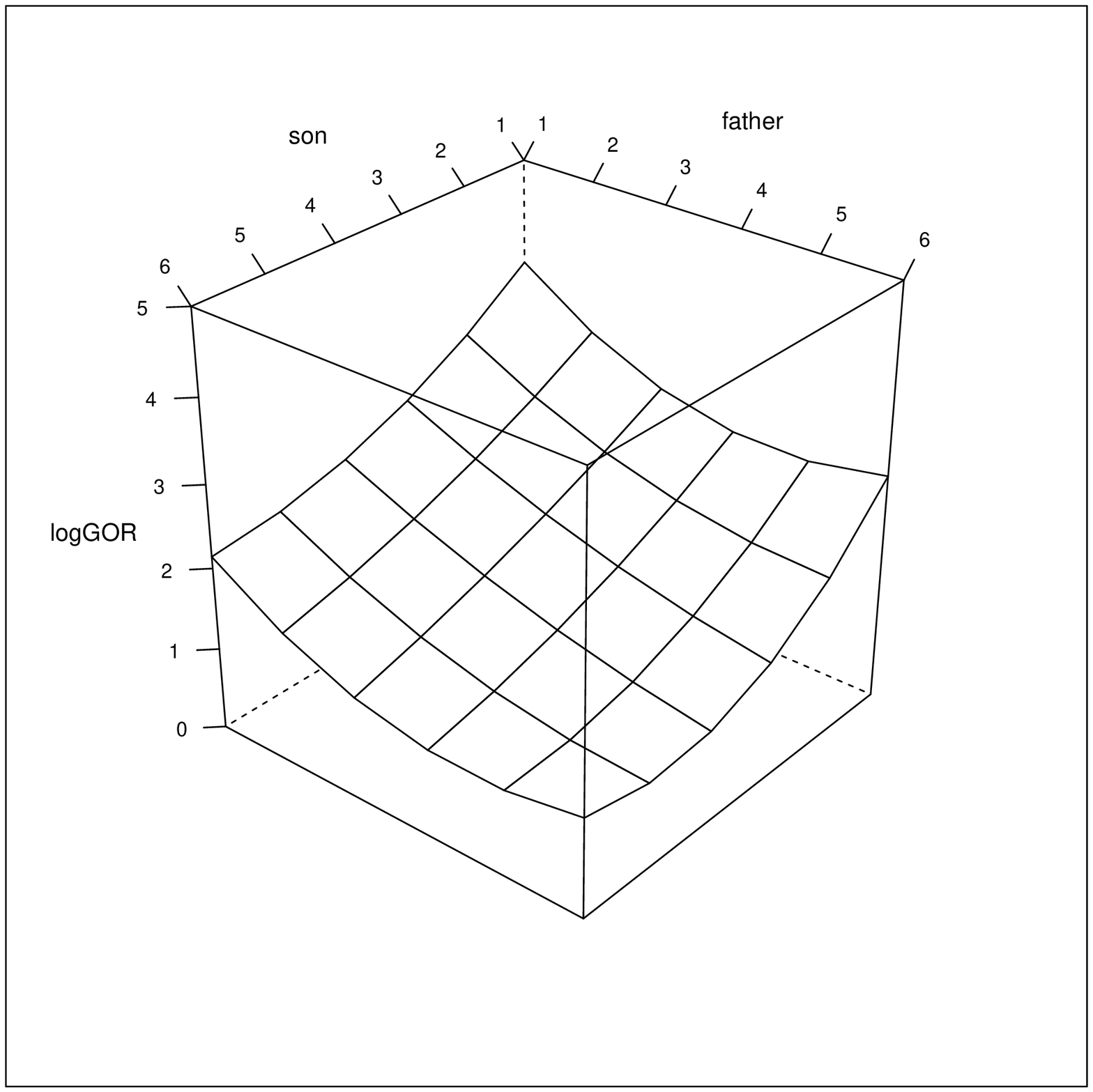}\hspace{-.15cm}
\includegraphics[width=7cm, height=7.5cm]{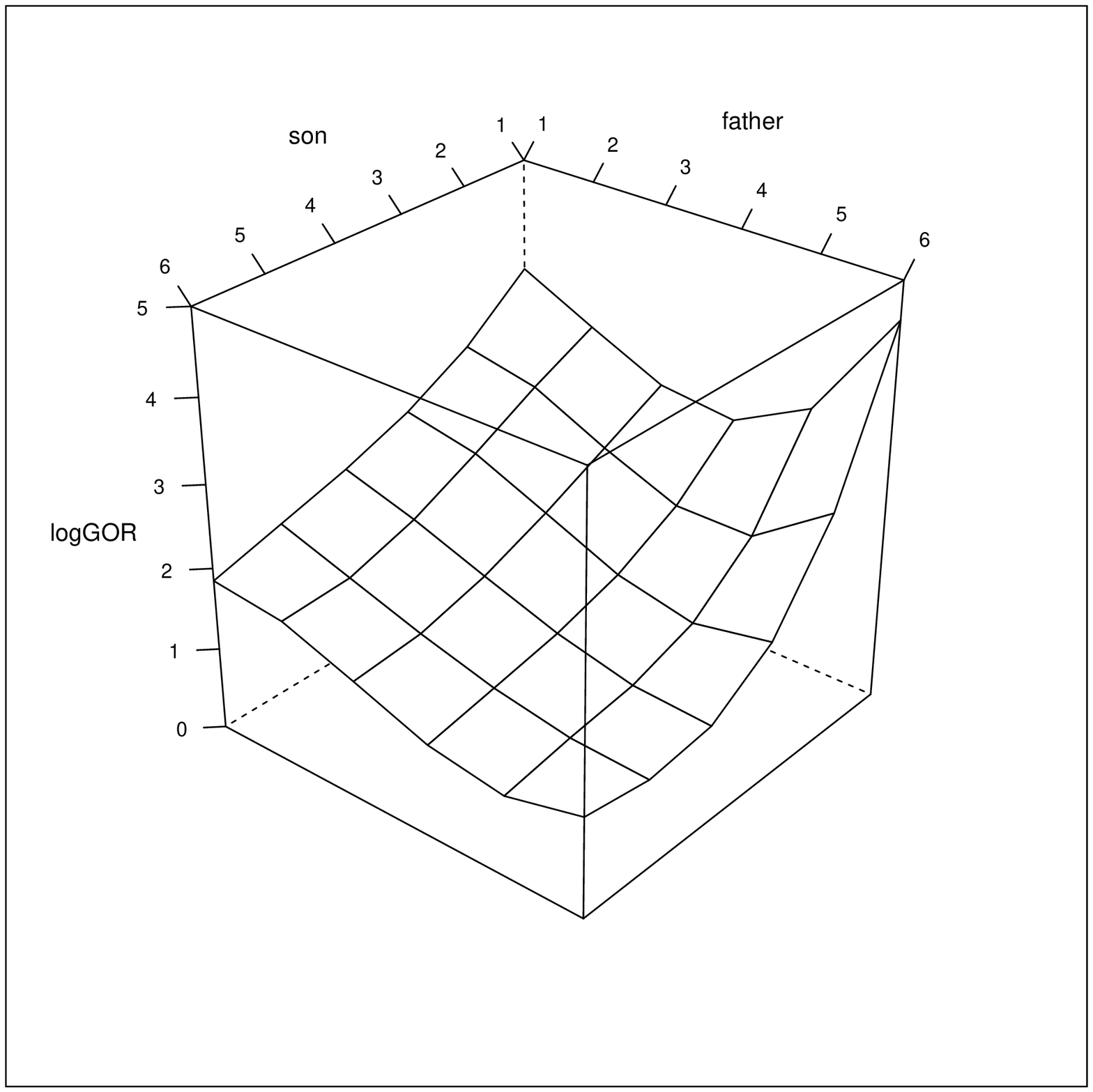}\vspace{-.5cm}\\
\caption{Association structures for the OS data set. In lexicographical order: the observed (top-left) and fitted log-GORs by the interaction of polynomials of second degree with non-integer scores (top-right), second (bottom-left) and third (bottom-right)degree with integer scores.}\label{fig:7}
\end{figure*}
\noindent The saturated model for the joint distribution involves 48 parameters: 6 global logits for each marginal and 36 log-GORs. Since the interest is in modeling the association structure, we concentrate only on the 36 log-GORs. Figure \ref{aic4} shows the AIC for the NUPOM for varying smoothing parameter and different orders of penalization.
\par Due to the symmetry of the association structure, the difference operator orders $s_{3}$ and $s_{4}$ are assumed to be equal and indicated by $s$. The model with first order penalization, whose $AIC$ tends to the UPOM value, is clearly inadequate. The minimum value of $AIC$ is $22236.65$, corresponding to $\log(\lambda)=4$ and $s=3$. This represents a model with a smoothed association structure which tends to a surface defined by row and column interactions of second degree polynomials. Thus, on the grounds of $AIC$ only, one could choose this model. However, observe that for high values of $\lambda$, the models with $s=3\,(AIC=22239.96)$ and $s=4\,(AIC=22238.34)$ provide good fits as well, with a slight evidence in favor of the latter model, corresponding to a surface of third degree polynomials. The AIC differences of such models, with respect to the minimum AIC, are respectively 3.31 and 1.69, which are quite small \citep[p. 48]{Burnham:2000}. When it is possible, as in this case, it is preferable to choose a model providing integer and equally-spaced scores, on the grounds of greater interpretability and for the possibility to use classical test statistics whose asymptotic null distributions are well-known. Therefore, we report in Table \ref{mydeviance} the results (in terms of $AIC$ and deviance $G^2$) for the four polynomial models evaluated at the largest value of $\log(\lambda)=15$, along with the independence and saturated models.

\begin{table}[H]
\small
%\vspace{-.8cm}
\caption{Model selection based on AIC and deviance $G^2$ for the OS data set based on a BOLM with penalty term ARC2. The asterisks indicate a non significant difference with the saturated model.}
%\vspace{-.1cm}
\begin{center}
\begin{tabular}{llrrr}
  \hline
Model & Description & df & AIC & $G^2$\\
  \hline
1 & Independence & 36 & 23081.12 & 897.52\\
2 & Uniform association & 35 & 22392.83 & 207.22\\
3 & First degree polynomials & 32 & 22247.46 &  55.85\\
4 & Second degree polynomials & 27   & 22239.96 &*38.36\\
5 & Third degree polynomials & 20   & 22238.34 &**22.74\\
6 & Saturated & 0   & 22255.60 &   0.00 \\
   \hline
\end{tabular}
\end{center}
\vspace{-.4cm}
\begin{flushleft}
\hspace{.2cm}
*p=0.07, **p=0.3.
\end{flushleft}
\label{mydeviance}
\end{table}

Model 1 has been fitted by using the ridge-type penalty term, such that the estimated global log-odds ratios tend to zero for high values of smoothing parameter. Model 4 provides the most parsimonious but yet acceptable fit, with only 9 estimated parameters and p-value = 0.07. Model 5 estimates only 16 parameters, providing a comparable fit ($G^2=22.74$), with a not significant difference with the saturated model (p-value = 0.3). This means the ordinal association structure of occupational status can be well fitted by two polynomials of second or third degree, that is a polynomial surface of fourth or sixth degree, respectively. Notice that Models 4 and 5 are more parsimonious than the best model found in Lapp \emph{et al.}'s analysis, i.e. the Dale model, including row effects, column effects, and interactions, while maintaining a comparable fit in terms of $G^2$. The observed structure of global log-odds ratios and the selected polynomial models are graphically showed in Figure \ref{fig:7}.\par  As we can see, the ordinal association structure is always positive, but it decreases as both the social statuses increase. Observe that the second degree polynomials using integer and non-integer scores show very slight differences. Finally, notice that also Lapp \emph{et al.} hypothesized the possibility to fit a ``symmetric second degree polynomial'' model.

\subsection{The liver disease patients data set}\label{liver}
The data set consists of 256 records related to the liver disease progression of patients. The two outcomes, both measured
on the same day, are the liver biopsy (named \emph{STAGE}), considered the natural gold standard, and a categorized version of transient elastography (\emph{STIFF}) according to cutoffs suggested by \cite{Castera:2005}, to measure liver stiffness. Both the responses have three ordered categories, \emph{STIFF} with levels 1,2 and 3, which correspond to stiffness \\classes $[0,7.1)$, $[7.1,12.5)$ and $[12.5,\infty)$, respectively; as for \emph{STAGE},the initial five categories F0-F4 have been collapsed as follows: 1, corresponding to the biopsy stage $<F2$, 2 ($=\{F2,F3\}$) and 3 ($=F4$).
Aim of the study is to evaluate the concordance between the outcomes in order to find profiles of ``discordant'' patients. This is done using the bivariate logistic model, by employing the log global odds as marginal parameters and the log global odds ratio as association measure. For these data, a first analysis with dichotomized responses was made by \cite{Calvaruso:2009} but estimation problems inhibited an analysis with 3 levels for each response. Table \ref{test} shows the cross-classification of the responses, ignoring covariates.\\

\begin{table}[htbp]
\small
\caption{Marginal Cross-classification of the responses and empirical global log-odds ratios for the liver disease patient data.}
\centering
\begin{tabular}{lccc}
& \multicolumn{3}{c}{\emph{STIFF}}\\
%\cline{2-4}
\emph{STAGE} & 1 & 2 & 3\\
\hline
1 & 71 & 20 & 0\\
  & (1.72) & ($+\infty$) &\\
  &&&\\
2 & 56 & 20 & 8\\
  & (3.18) & (3.31) &\\
  &&&\\
3 & 8 & 27 & 46\\
\hline
\end{tabular}
\label{test}
\end{table}

\par From a first look at Table \ref{test}, it is possible to notice an overall positive association between the responses, even if there are many discordant patients, mainly the fifty-six in $STAGE=2,\,STIFF=1$. Among the covariates, the patient's sex (\emph{SEX}), age (\emph{AGE}), alanine aminotransferase  (\emph{ALT}) measured in $U/L$ and platelet (\emph{PLT}) levels measured in $10^3$ mmc,  are considered. For modelling purposes, the covariates were centered with respect to their means and, after a backward selection, the following sets of variables have been considered:\\

$$\mathcal{P}_{1}^{0}=\{\textrm{Intercept}_{\overline{123}},SEX_{123},\,AGE_{123},\,ALT_{123},\,PLT_{1\overline{2}3}\},$$
$$\mathcal{P}_{2}^{0}=\{\textrm{Intercept}_{\overline{123}},SEX_{123},\,AGE_{123},\,ALT_{123},\,PLT_{123}\},$$
$$\mathcal{P}_{3}^{0}=\{\textrm{Intercept}_{\overline{12}3},SEX_{123},\,AGE_{123},\,ALT_{123},\,PLT_{1\overline{2}3}\},$$
$$\mathcal{P}_{4}^{0}=\{\textrm{Intercept}_{\overline{12}3},SEX_{123},\,AGE_{123},\,ALT_{123},\,PLT_{123}\},$$
$$\mathcal{P}_{5}^{0}=\{\textrm{Intercept}_{\overline{123}},AGE_{123},\,ALT_{123},\,PLT_{1\overline{2}3}\},$$
$$\mathcal{P}_{6}^{0}=\{\textrm{Intercept}_{\overline{123}},AGE_{12},\,ALT_{123},\,PLT_{1\overline{2}}\},$$\\

\noindent where $\textrm{Intercept}_{\overline{123}}$ indicates that the marginal and association intercepts are category dependent, whereas $\textrm{Intercept}_{\overline{12}3}$ indicates that the intercepts for the association are category independent. Here, to indicate that variable $AGE$ is included both in marginal and association predictors, we use $AGE_{123}$, whereas $AGE_{12}$ indicates that such variable is included only in the marginal predictors.
\par Computational problems have arisen when we tried  to estimate the models which assumed a non-uniform association. We have overcome such problems by regularizing the parameter space of the association intercepts. In particular, we have employed the ARC1 penalty term, with smoothing parameter $\lambda_3=1$ which is the minimum common penalization value for which Fisher scoring did not fail in all the fitted models of $NU(P)POM$ type. By considering the results from the simulation in Section \ref{section:testing}, we decided to use a $\chi^2$ distribution to approximate the $LR_P$ asymptotic distribution. This is done in the model selection reported in Table \ref{tabaic4}.

\begin{table}[ht]
\setlength{\tabcolsep}{4pt}
\small
\caption{Model selection based on the $AIC$ and the $LR_P$ statistic for the liver disease patients data using the ARC1 penalty term.}
\begin{center}
\begin{tabular}{llrcrrrr}

  \hline
Model & Description & \# par. & AIC & vs  & $LR_P$  & df & p-value\\
  \hline
1&$NUPPOM(\mathcal{P}_{1}^0)$ & 21 & 879.5 & - & -     & - & -\\
2&$NUPOM(\mathcal{P}_{2}^0)$  & 20 & 882.1 & 1 & 4.67  & 1 & $0.031$ \\
3&$UPPOM(\mathcal{P}_{3}^0)$  & 18 & 895.7 & 1 & 16.51 & 3 & $<0.001$ \\
4&$UPOM(\mathcal{P}_{4}^0)$   & 17 & 898.2 & 1 & 20.99 & 4 & $<0.001$ \\
\multicolumn{8}{c}{\hdashrule[0.5ex]{8cm}{.5pt}{.8pt}}\\
5&$NUPPOM(\mathcal{P}_{5}^0)$ & 18 & 874.7 & 1 & 1.35  & 3 & $0.719$ \\
6&$NUPPOM(\mathcal{P}_{6}^0)$ & 17 & 873.0 & 5 & 0.36  & 1 & $0.549$ \\
7&$NUPPOM(\mathcal{P}_{7}^0)$ & 16 & 872.7 & 6 & 1.23  & 1 & $0.268$ \\
7&- & - & - & 1 & 2.93  & 5 & $0.711$ \\
   \hline
\end{tabular}
\end{center}
\label{tabaic4}
\end{table}

\par For each model we have selected, the table reports its description, the number of estimated parameters and the AIC. The next columns refer to comparisons between nested models, specified by the column headed \lq\lq vs". The last three columns report the results of such a comparison in terms of penalized log-likelihood ratio statistic, along with degrees of freedom and p-values. Before proceeding to variable selection we first checked, for each variable, the hypothesis $UPOM$, versus several alternatives: $UPPOM$, $NUPOM$, $NUPPOM$. The table reports such comparisons for Models 1-4. Model 1 is the most complex model we have considered, a $NUPPOM$ defined on set $\mathcal{P}_{1}^{0}$. This model assumes that the effect of variable $PLT$ on $STIFF$ depends on the categories of such response variable. Models 2-4 represent hypotheses of uniform association and/or (partially) proportional odds, and these models are compared to Model 1, for which none of these hypotheses holds. Although the $LR_P$ test for model comparison is approximated, some results seem to be clear. For example, the difference between model 1 and 3 (or 4) is highly significant, so we can claim that the hypothesis of $UPPOM$ (or $UPOM$) does not hold. Models 5-7 concern a backward model selection starting from Model 1.  The last row reports the comparison between Models 7 and 1, for which the difference  between the starting model and the final model is not significant (p-value=0.711). In model 7, variable $ALT$ is the only one which has significant (global) effect for the association model. By AIC, the model with the best trade-off between goodness-of-fit and parsimony is still Model 7. Estimates for this final model are reported in Table \ref{tabliver}.
\par For both outcomes, variables $ALT$, $AGE$ and $PLT$ are significant. In particular the platelet level has a category-dependent effect for $STIFF$ which is higher for the log global odds 1-2 than 3. In particular, a patient at older age, higher $ALT$ and lower $PLT$ values is more at risk of having a greater liver stiffness than a patient with mean values. In addition, the $ALT$ effect for $STIFF$ is about twice as strong as for biopsy stage. The effect of $ALT$ in the association is significant, and considering the intercepts values as well, higher $ALT$ values imply a global reduction of the association, especially for individuals in class $STIFF<7.1$ and $STAGE=1$.

\begin{table}[H]
 \setlength{\tabcolsep}{5pt}
\small
\caption{Estimates for Model 4.}
\begin{center}
\begin{tabular}{llrrrr}
  \hline
   response &variable & estimate & se & z & p.value \\
  \hline
\\
STAGE & Intercept 1 & -0.7398 & 0.1417 & -5.2217 & 0.0000 \\
      & Intercept 2 & 0.8798 & 0.1452 & 6.0610 & 0.0000 \\
      & ALT & -0.0047 & 0.0018 & -2.6217 & 0.0088 \\
      & AGE & -0.0407 & 0.0105 & -3.8985 & 0.0001 \\
      & PLT & 0.0093 & 0.0021 & 4.4209 & 0.0000 \\
&&&&&\\
STIFF & Intercept 1 & 0.1021 & 0.1377 & 0.7421 & 0.4581 \\
      & Intercept 2 & 1.7861 & 0.2003 & 8.9152 & 0.0000 \\
      & ALT & -0.0089 & 0.0019 & -4.6721 & 0.0000 \\
      & AGE & -0.0422 & 0.0110 & -3.8274 & 0.0001 \\
      & PLT 1 & 0.0087 & 0.0024 & 3.6768 & 0.0002 \\
      & PLT 2 & 0.0154 & 0.0030 & 5.1944 & 0.0000 \\
&&&&&\\
STAGE&&&&&\\
vs STIFF & Intercept & 1.5164 & 0.3320 & 4.5680 & 0.0000 \\
      & Intercept & 3.3436 & 0.6207 & 5.3869 & 0.0000 \\
      & Intercept & 2.8285 & 0.3880 & 7.2892 & 0.0000 \\
      & Intercept & 2.8134 & 0.4189 & 6.7161 & 0.0000 \\
      & ALT & -0.0075 & 0.0036 & -2.0668 & 0.0387 \\
 \\
\hline
\end{tabular}
\end{center}
\label{tabliver}
\end{table}

\section{Discussion}\label{sec5}
\noindent We have shown how to fit a BOLM by penalized ML estimation with some penalty terms for a ``vertical penalization'', that is across response levels. Particular emphasis on the terms ARC1 and ARC2, penalizing adjacent row and column effects, has been given. The motivation for such an approach is, on one hand, its flexibility in modeling situations in which ML estimation by traditional Fisher scoring appears somewhat difficult, and on the other hand, the possibility to consider the fit of a $NUPPOM$, which lies between a $UPOM$, which may give a poor fit, and a $NUNPOM$, often less useful and somewhat more complicated to estimate than a $UPOM$. The penalized log-likelihood ratio $LR_P$ statistic has been considered to check the hypothesis that certain effects are category independent.  To our knowledge, the asymptotic distribution of $LR_P$ for the considered hypothesis is not known, though we have shown, by simulation, that for relatively small smoothing values the $\chi^2$ may be a good approximation. However, as far as the distributional properties of penalized likelihood ratio test-statistics are concerned, further investigations are necessary. The potential of penalized estimates by penalty term ARC1 has been shown by simulation and by an application to an original data set. In addition, the BOLM has been fitted using the penalty term ARC2 to a literature data set  for comparison with the alternative Dale and Goodman RC models, showing parsimony  while preserving a satisfactory the goodness-of-fit. In some sense,  ARC2 generalizes ARC1, permitting to fit restricted versions of the Dale model, by inserting row or column effects, but also polynomial effects models, with scores chosen by data. All the computations, including the model implementation, have been carried out by an original R code which can be requested from the authors. Furthermore, an R package implementing all the methods proposed in this paper is under preparation.
\section*{Appendix A: penalized maximum likelihood estimation}
\noindent Let $\partial l/\partial\bm{\pi}_{i}=\textrm{diag}(\bm{\pi}_{i})^{-1}\bm{y}_{i}$,
$\partial \bm{\pi}_{i}/\partial\bm{\eta}_{i}=(\mathbf{C}'\mathbf{D}_{i}^{-1}\mathbf{L})^{-1}$ and
 $\partial \bm{\eta}_{i}/\partial\bm{\beta}=\mathbf{X}_{i}$. By using the chain rule, the first derivative of the penalized log likelihood with respect to $\bm{\beta}$ is
 $$\frac{\partial l_{P}}{\partial\bm{\beta}}=\sum_{i=1}^{m}{\frac{\partial l}{\partial\bm{\pi}_{i}}
\frac{\partial\bm{\pi}_{i}}{\partial\bm{\eta}_{i}}\frac{\partial\bm{\eta}_{i}}{\partial\bm{\beta}}}-\mathbf{P}\bm{\beta},$$
the \emph{penalized score function} is
$$\bm{s}_{P}(\bm{\beta};\bm{y}_{i})=\sum_{i=1}^m{[(\mathbf{C}'\mathbf{D}_{i}^{-1}\mathbf{L})^{-1}\mathbf{X}_{i}]'
\textrm{diag}(\bm{\pi}_{i})^{-1}\bm{y}_{i}}-\mathbf{P}\bm{\beta},$$
and the \emph{penalized Fisher matrix} is  $$\mathbf{F}_{\hspace{-.02cm}P}(\bm{\beta})\hspace{-.03cm}=\hspace{-.03cm}\sum_{i=1}^{m}{n_{i}\mathbf{X}'_{i}(\mathbf{L}'\mathbf{D}_{i}'^{-1}\mathbf{C})^{-1}
\textrm{diag}(\bm{\pi}_{i})^{-1}(\mathbf{C}'\mathbf{D}_{i}^{-1}\mathbf{L})^{-1}\mathbf{X}_{i}+\mathbf{P}}.$$

\noindent Using these formulas, the $(k+1)$th iteration of the Fisher scoring is
$\hat{\bm{\beta}}^{(k+1)}=\hat{\bm{\beta}}^{(k)}+l_{step}\mathbf{F}_{P}(\hat{\bm{\beta}}^{(k)})^{-1}\bm{s}_{P}(\hat{\bm{\beta}}^{(k)})$, where $l_{step}$ is a positive scalar representing the step length. Since the iterative procedure may produce incompatible $\bm{\beta}$ values for $\bm{\pi}$, a value smaller than 1 for $l_{step}$, say $0.5$ or smaller, may be necessary, even if this inevitably increases the number of iterations. As a reasonable starting value for $\bm{\beta}$, one could set to zero the regression coefficients corresponding to covariates, together with the global log-odds ratios intercepts, whereas the global logits intercepts have to be chosen by taking into account the inequality constraints (\ref{constr}). The variance\ covariance matrix of $\hat{\bm{\beta}}$ is given by $V(\hat{\bm{\beta}})\hspace{-.03cm}=\hspace{-.03cm}F_{\hspace{-.02cm}P}(\hat{\bm{\beta}})^{-1}$. When a $NUPPOM$ is considered, the form of matrix $\mathbf{X}_{i}$ is $\mathbf{X}_{i}=\bigoplus_{k=1}^{3}{\mathbf{X}_{k,i}}$, where:
\begin{displaymath}
\mathbf{X}_{k,i}=\left(
\begin{array}{ccccccccc}
1 &        & 0      & \bm{x}'_{i,\mathcal{S}_{k}} & \bm{x}'_{i,\bar{\mathcal{S}}_{k}}&   & \bm{0}_{}'\\
  & \ddots &        & \vdots                     &                                 & \ddots &\\
0 &        & 1      & \bm{x}'_{i,\mathcal{S}_{k}} &  \bm{0}' & &\bm{x}'_{i,\bar{\mathcal{S}}_{k}}\\
\end{array}
\right).
\end{displaymath}
\noindent Thus the full design matrix is simply $\mathbf{X}=(\mathbf{X}_{1}',\mathbf{X}_{2}',...,\mathbf{X}_{m}')'$. The \emph{weight function} for the $i$th observation is defined as $\mathbf{W}_{i}(\bm{\beta})=n_{i} \left(\frac{\partial \bm{\pi}_{i}}{\partial\bm{\eta}_{i}'}\right)\textrm{diag}(\bm{\pi}_{i})^{-1}\left(\frac{\partial \bm{\pi}_{i}}{\partial\bm{\eta}_{i}}\right)$, the weight matrix is
$\mathbf{W}(\bm{\beta})=(\mathbf{W}_1(\bm{\beta})',\mathbf{W}_2(\bm{\beta})',\dots, \mathbf{W}_m(\bm{\beta})')'$, the \emph{hat matrix} is $\mathbf{H}=\mathbf{X}(\mathbf{X}'\mathbf{W}(\hat{\bm{\beta}})\mathbf{X}+\mathbf{P})^{-1}\mathbf{X}'\mathbf{W}(\hat{\bm{\beta}})$, and the \emph{Akaike Information Criterion} is $AIC=-2(l(\hat{\bm{\beta}})-tr(\mathbf{H}))$.

\section*{Appendix B: penalty terms in matrix form}
\label{ptmf}
When (\ref{sis3}) or (\ref{sis4}) is used
$\mathbf{P}=\mathbf{E}'\boldsymbol{\mathrm{\Lambda}}'^{1/2}\boldsymbol{\mathrm{\Lambda}}^{1/2}\mathbf{E},$
where $\boldsymbol{\mathrm{\Lambda}}$ is the matrix of smoothing values, and $\mathbf{E}=\bigoplus_{k=1}^3\mathbf{E}_{k}$.\\ \par In (\ref{sis3}) matrices  $\boldsymbol{\mathrm{\Lambda}}$ and $\mathbf{E}_{k}$ depend on the penalty (ridge or ARC1). Let $\bm{d}=(D_1-1,D_2-1,(D_1-1)(D_2-1))'$ a vector indexed by $d_{k},\,k=1,2,3,$ and let $\bm{\lambda}_{k_{|\bar{\mathcal{S}}_{k}^0|}}'=(\lambda_{k,0},\bm{\lambda}'_{k_{|\bar{\mathcal{S}}_{k}|}}) $ be the smoothing values vector of length $|\bar{\mathcal{S}}_{k}^0|,$ that is the cardinality of the set of variables undergone to penalization for the $k$th equation in (\ref{sis2}). Then $$\boldsymbol{\mathrm{\Lambda}}=\textrm{diag}(\bm{\lambda}_{1,\mathcal{P}_1^0}',\bm{\lambda}_{2,\mathcal{P}_2^0}', \bm{\lambda}_{3,\mathcal{P}_3^0}'),$$ where
$\bm{\lambda}_{k,\mathcal{P}^0_k}'=(\lambda_{k,0}\bm{1}'_{(d_k-1)},\bm{0}'_{_{|\mathcal{S}_{k}|}},
   \bm{\lambda}'_{k_{|\bar{\mathcal{S}}_{k}|}}\varotimes\bm{1}'_{(d_k-1)}).$
\subsection*{\emph{Ridge-type penalty term}}
\par Let $\mathbf{I}_{d_k}$ be the $d_k\times d_k$ identity matrix. Then
$$\mathbf{E}_k=(\mathbf{I}_{d_k},\bm{0}_{d_k\times|\mathcal{S}_k|},
\bm{1}'_{|\bar{\mathcal{S}}_k|}\varotimes \mathbf{I}_{d_k}).$$
\subsection*{\emph{ARC1 penalty term}}
\par For the ARC1 penalty, let $\mathbf{T}_{k}$ be the $d_k\times d_k$ upper triangular matrix of 1's. Its inverse $\mathbf{V}_{k}=\mathbf{T}_{k}^ {-1}$ has 1's on the main diagonal, -1's on the first superdiagonal and 0's elsewhere. Further, let $\mathbf{V}_{k_{-1}}$ be the matrix $\mathbf{V}_{k}$ ignoring the last row. Then $$\mathbf{E}_k=(\mathbf{V}_{k_{-1}},\bm{0}_{(d_k-1)\times|\mathcal{S}_k|},
\bm{1}'_{|\bar{\mathcal{S}}_k|}\varotimes \mathbf{V}_{k_{-1}}).$$
\subsection*{\emph{ARC2 penalty term}}
\par  Let $\bm{c}=(D_1-1,D_2-1,D_2-1,D_1-1)'$ a vector indexed by $c_{h},\,h=1,...,4$, and let $k=1,2,3$. Then
$$\boldsymbol{\mathrm{\Lambda}}=\textrm{diag}(\bm{\lambda}_{1,\mathcal{P}_1^0}',\bm{\lambda}_{2,\mathcal{P}_2^0}',
  \bm{\lambda}_{3,\mathcal{P}_3^0}',\bm{\lambda}_{4,\mathcal{P}_3^0}'),$$ where
$\bm{\lambda}_{h,\mathcal{P}_k^0}'=(\lambda_{h,0}\bm{1}'_{(c_h-1)},\bm{0}'_{_{|\mathcal{S}_{k}|}},
 \bm{1}'_{(c_h-1)}\varotimes \bm{\lambda}'_{h_{|\bar{\mathcal{S}}_{k}|}}).$
%and $\bm{\lambda}_{h_{|\bar{\mathcal{S}}_k^0|}}'=(\lambda_{h,0},$ $\bm{\lambda}'_{h_{|\bar{\mathcal{S}}_{k}|}}).$

Define $s_{h,j},\,\,j \in \bar{\mathcal{S}}^0,$ the order of operator $\Delta^{s_{h,j}}$,  for the $j$th variable, also including the intercepts.  Let $\mathbf{T}_{h}$ be the $c_h\times c_h$ upper triangular matrix of $1's$ and let $\mathbf{V}_{h}=\mathbf{T}_{h}^{-1}$. Let $\mathbf{V}^{s_{h,j}}=\prod_{h=1}^{s_{h,j}}-\mathbf{V}_{h}$,
let $\mathbf{V}^{s_{h,j}}_{-s_{h,j}}$ be the matrix $\mathbf{V}^{s_{h,j}}$ ignoring the last $s_{h,j}$ rows and let $\mathbf{U}_{k,\bar{\mathcal{S}}^0}=(\mathbf{U}_{k,0},\mathbf{U}_{k,\bar{\mathcal{S}}}),$ where
$\mathbf{U}_{k,\bar{\mathcal{S}}}=(\mathbf{U}_{k,1},...,\mathbf{U}_{k,|\bar{\mathcal{S}}|})$ and
$\mathbf{U}_{1,j}=\mathbf{V}^{s_{1,j}}_{-s_{1,j}}$,
$\mathbf{U}_{2,j}=\mathbf{V}^{s_{2,j}}_{-s_{2,j}}$ and
$\mathbf{U}_{3,j}=(\mathbf{I}_{(c_3)}\varotimes \mathbf{V}^{s_{3,j}}_{-s_{3,j}})\varoplus
(\mathbf{V}^{s_{4,j}}_{-s_{4,j}}\varotimes \mathbf{I}_{(c_4)}).$ Then   $$\mathbf{E}_{k}=(\mathbf{U}_{k,0},\bm{0}_{(d_k-1)\times|\mathcal{S}_k|},
\mathbf{U}_{k,\bar{\mathcal{S}}}).$$

\subsection*{\emph{The penalty term for ordering constraints}}
\par For (\ref{sis_constraints}) let $\mathbf{N}=(\mathbf{1}_{n}\varotimes \mathbf{E})'$, being $\mathbf{E}$ defined as for ARC1, and $n$ the sample size. Let $\boldsymbol{\mathrm{\Lambda}}=(\lambda_{1}\mathbf{I}_{n},\lambda_{2}\mathbf{I}_{n},\bm{0}_{n \times n})$.  Then $$\mathbf{P}=\mathbf{X}'\mathbf{N}'\boldsymbol{\mathrm{\Lambda}}'^{1/2}I(\bm{\beta}'\mathbf{X}'\mathbf{N}')
I(\mathbf{N}\mathbf{X}\bm{\beta})\boldsymbol{\mathrm{\Lambda}}^{1/2}\mathbf{NX},$$ where $I(\mathbf{NX}\bm{\beta}\leq 0)$ is element-wise, that is $I(a_{ij}\leq 0)=1$ if true, 0  otherwise.

\begin{acknowledgements}
This paper has been supported by the Italian Ministerial grant
PRIN 2008 ``Measures, statistical models and indicators for the
assessment of the University System'', n. 2008WXSLH.
\end{acknowledgements}

\bibliographystyle{spbasic}      % basic style, author-year citations
\bibliography{bibtesi-1abbr}   % name your BibTeX data base

\end{document}